\input ./style/arxiv-general.cfg
\documentclass[aos,MSNbibl,nameyear,dvips]{arximspdf}
\makeatletter
   \@ifpackageloaded{graphicx}{}{\usepackage{graphicx}}
\makeatother
\usepackage{url,breakurl}

\doi{10.1214/15-AOS1344}
\volume{43}
\issue{5}
\pubyear{2015}
\firstpage{2102}
\lastpage{2131}
\docsubty{FLA}

\makeatletter

\newcommand{\rrvert}{\vert}

\newcommand{\rrVert}{\Vert}
\newcommand{\llvert}{\vert}
\newcommand{\llVert}{\Vert}
\newcommand{\eqref}[1]{(\ref{#1})}

\newtheorem{theorem}{Theorem}

\newtheorem{proposition}{Proposition}

\newtheorem{lemmaa}{Lemma}
\newproclaim{example}{Example}
\newproclaim{assumption}{Assumption}
\newproclaim{remark}{Remark}
\makeatother

\begin{document}
\begin{frontmatter}

\title{Estimation and inference in generalized additive coefficient
models for nonlinear interactions with high-dimensional covariates}
\runtitle{GACM for nonlinear interactions}

\begin{aug}
\author[A]{\fnms{Shujie}~\snm{Ma}\thanksref{T1,m1}\ead[label=e1]{shujie.ma@ucr.edu}},
\author[B]{\fnms{Raymond J.}~\snm{Carroll}\thanksref{T2,m2,mk}\ead[label=e2]{carroll@stat.tamu.edu}},
\author[C]{\fnms{Hua}~\snm{Liang}\corref{}\thanksref{T3,m3}\ead[label=e3]{hliang@gwu.edu}}
\and
\author[D]{\fnms{Shizhong}~\snm{Xu}\thanksref{m1}\ead[label=e4]{shizhong.xu@ucr.edu}}
\runauthor{Ma, Carroll, Liang and Xu}
\thankstext{T1}{Supported in part by NSF Grant DMS-13-06972.}
\thankstext{T2}{Supported by  National Cancer Institute Grant
U01-CA057030.}
\thankstext{T3}{Supported in part by NSF Grants DMS-14-40121 and
DMS-14-18042 and by Award Number
11228103, made by National Natural Science Foundation of China.}

\affiliation{University of California, Riverside\thanksmark{m1},
Texas A\&M University\thanksmark{m2}, University of Technology
Sydney\thanksmark{mk} and
George Washington University\thanksmark{m3}}
\address[A]{S. Ma\\
Department of Statistics\\
University of California, Riverside\\
Riverside, California 92521\\
USA\\
\printead{e1}}
\address[B]{R. J. Carroll\\
Department of Statistics\\
Texas A\&M University\\
College Station, Texas 77843\\
USA\\
and\\
School of Mathematical Sciences\hspace*{8pt}\\
University of Technology\\
Sydney, Broadway NSW 2007\\
Australia\\
\printead{e2}}
\address[C]{H. Liang\\
Department of Statistics\\
George Washington University\\
Washington, DC 20052\\
USA\\
\printead{e3}}
\address[D]{S. Xu\\
Center for Plant Cell Biology\\
University of California, Riverside\\
Riverside, CA 92521\\
USA\\
\printead{e4}}
\end{aug}

%
\received{\smonth{9} \syear{2014}}
%
\revised{\smonth{5} \syear{2015}}

%
\begin{abstract}
In the low-dimensional case, the generalized additive coefficient model
(GACM) proposed by Xue and Yang [\textit{Statist. Sinica} \textbf{16}
(2006) 1423--1446]
has been demonstrated to be a
powerful tool for studying nonlinear interaction effects of variables.
In this paper, we propose estimation and inference procedures for the GACM
when the dimension of the variables is high. Specifically, we propose a
groupwise penalization based procedure to distinguish significant
covariates for the ``large $p$ small $n$'' setting. The procedure is
shown to be consistent for model structure identification. Further, we
construct simultaneous confidence bands for the coefficient functions
in the selected model based on a refined two-step spline estimator. We also
discuss how to choose the tuning parameters. To estimate the standard
deviation of the functional estimator, we adopt the smoothed bootstrap
method. We conduct simulation experiments to evaluate the numerical
performance of the proposed methods and analyze an obesity data set
from a
genome-wide association study as an illustration.
\end{abstract}


%
\end{frontmatter}

\section{Introduction}\label{sec:intr}

Regression analysis is a commonly used statistical tool for modeling
the relationship between a scalar dependent variable $Y$ and one or
more explanatory variables denoted as {$\mathbf{T}=(T_{1},T_{2},\ldots,T_{p})
^{\mathrm{T}}$.} To study the marginal effects of the
predictors on the response, one may fit a generalized linear model (GLM),
%
\begin{equation}
E(Y|\mathbf{T})=\mu(\mathbf{T})=g^{-1}\bigl\{\eta(\mathbf{T})\bigr\},\qquad
\eta (\mathbf{T}%
)=\sum_{\ell=1}^{p}
\alpha_{\ell0}T_{\ell}, \label{MOD:GLM}
\end{equation}
where $g$ is a known monotone link function, and $\alpha_{\ell0}$,
$1\leq
\ell\leq p$, are unknown parameters. Sometimes, the effect of one variable
may change with other variables; that is, there is an interaction
effect. { By
letting $T_{1}=1$, to incorporate} the interaction effects of $\mathbf
{T}$ and the other variables,
denoted as $\mathbf{X}=(X_{1},\ldots,X_{d})^{\mathrm{T}}$, model (\ref
{MOD:GLM}) can be modified to $E(Y|\mathbf{X},\mathbf{T})=\mu(\mathbf{X},\mathbf{T}%
)=g^{-1}\{\eta(\mathbf{X},\mathbf{T})\}$ with
%
\begin{equation}
\eta(\mathbf{X},\mathbf{T})=\alpha_{10}+\sum
_{\ell=2}^{p}\alpha _{\ell
0}T_{\ell}+
\sum_{k=1}^{d}\alpha_{1k}X_{k}+
\sum_{\ell
=2}^{p}\sum
_{k=1}^{d}\alpha_{\ell k}X_{k}T_{\ell},
\label{MOD:Linear2}
\end{equation}
where $\alpha_{\ell k}$ for $0\leq k\leq d$ and $1\leq
\ell\leq p$ are parameters. After a direct reformulation, model (\ref%
{MOD:Linear2}) can be written as
%
\begin{equation}
\eta(\mathbf{X},\mathbf{T})=\sum_{\ell=1}^{p}
\Biggl(\alpha_{\ell
0}+\sum_{k=1}^{d}
\alpha_{\ell k}X_{k}\Biggr)T_{\ell}. \label{MOD:Linear}
\end{equation}

\begin{figure}

\includegraphics{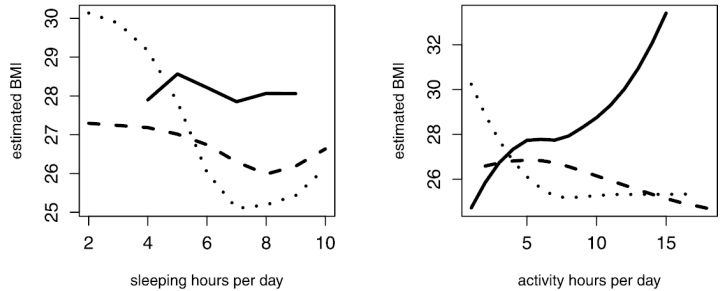}

\caption{Plots of the estimated BMI against sleeping
hours per day (left panel) and activity hours per day (right panel) for the
three genotypes AA (solid line), Aa (dashed line) and aa (dotted line) of
SNP rs242263 in the Framingham study, where A is the minor allele.}
\label{FIG:BMI}
\end{figure}

Here the effect of each $T_{\ell}$ changes linearly with $X_{k}$. However,
in practice, this simple linear relationship may not reflect the true
changing patterns of the coefficient with other covariates. We here use an
example of gene and environment (G${}\times{}$E) interactions for illustration.
It has been noticed in the literature that obesity is linked to genetic
factors. Their effects, however, can be altered under different
environmental factors such as sleeping hours [\citet{knutson:2012}] and
physical activity [\citet{Wareham.van.Ekelund:2005}]. To have a rough
idea of how
the effects of the genetic factors change with the environment, we explore
data from the Framingham Heart Study [\citet
{Dawber.Meadors.Moore:1951}]. In
Figure~\ref{FIG:BMI} we plot the estimated mean body mass index (BMI)
against sleeping hours per day and activity hours per day,
respectively, for
people with three possible genotype categories represented by AA, Aa
and aa,
and for one single nucleotide polymorphism (SNP). A detailed
description and
the analysis of this data set are given in Section~\ref{sec8}. We define
allele A as the minor (less frequent) allele. This figure clearly shows
different nonlinear curves for the three groups in each of the two
plots. By
letting $T_{\ell}$ be the indicator for the group $\ell$, the linear
function in model (\ref{MOD:Linear}) is clearly misspecified.

To relax the linearity assumption, we allow each $\alpha_{\ell k}X_{k}$
term to be an unknown nonlinear function of $X_{k}$, and thus extend
model (%
\ref{MOD:Linear}) to the generalized additive coefficient model (GACM)
%
\begin{equation}
\eta(\mathbf{X},\mathbf{T})=\sum_{\ell=1}^{p}
\Biggl\{\alpha_{\ell
0}+\sum_{k=1}^{d}
\alpha_{\ell k}(X_{k})\Biggr\}T_{\ell
}=\sum
_{\ell=1}^{p}\alpha_{\ell}(
\mathbf{X})T_{\ell}. \label{MOD:ADDITIVE}
\end{equation}
For identifiability, the functional components satisfy $E\{\alpha
_{\ell
k}(X_{k})\}=0$ for $1\leq k\leq d$ and $1\leq\ell\leq p$. The conditional
variance of $Y$ is modeled as a function of the mean, that is,
$\operatorname{var}(Y|\mathbf{X},\mathbf{T})=V\{\mu(\mathbf{X},\mathbf{T})\}=\sigma^{2}(\mathbf{X},\mathbf{T})$. In each
coefficient function of the GACM, covariates $X_{k}$ are continuous
variables. If some of them are discrete, they will enter linearly. For
example, if $X_{k}$ is binary, we let $\alpha_{\ell k}(X_{k})=\alpha
_{\ell
k}X_{k}$. In such a case, model (\ref{MOD:ADDITIVE}) turns out to be a
partially linear additive coefficient model. The linearity of (\ref%
{MOD:ADDITIVE}) in $T_{\ell}$ is particularly appropriate when those
factors are discrete, for example, SNPs in a genome-wide association
study (GWAS),
as in the data example of Section~\ref{sec8}.

For the low-dimensional case that the dimensions of $\mathbf{X}$ and $%
\mathbf{T}$ are fixed, estimation of model (\ref{MOD:ADDITIVE}) has been
studied; see \citet{xue.yang:2006,xue.liang:2010,liu.yang:2010} for a
spline estimation procedure and \citet{lee.mammen.park:2012} for a
backfitting algorithm. In modern data applications, model~(\ref
{MOD:ADDITIVE}%
), however, is particularly useful when $p$ is large. For example, in GWAS,
the number of SNPs, which is $p$, can be very large, but the dimension
of $
\mathbf{X}$ such as the environmental factors, which is $d$, is inevitably
relatively small. Moreover, the number of variables in $\mathbf{T}$ which
have nonzero effects is small. It therefore, poses new challenges to apply
model (\ref{MOD:ADDITIVE}) to the high-dimensional case including: (i)
how to
identify those important variables in $\mathbf{T}$, (ii) how to
estimate the
coefficient functions for the important covariates and (iii) how to conduct
inferences for the nonzero coefficient functions. For example, it is of
interest to know whether they are a function of a specific parametric form
such as constant, linear or quadratic, etc.

In the high-dimensional data setting, studying nonlinear interaction effects
has found much attention in recent years, and a few strategies have been
proposed. For example, \citet{jiang.liu:2014} proposed to detect variables
under the general index model, which enables the study of high-order
interactions among components of continuous predictors, which are
assumed to
have a multivariate normal distribution. Moreover, \citet{lian:2012}
considered variable selection in varying coefficient models which
allows the
coefficient functions to depend on one index variable, such as a
time-dependent variable.

When we would like to see how the effect of each genetic factor changes
under the influence of multiple environmental variables, the proposed
high-dimensional GACM (\ref{MOD:ADDITIVE}) becomes a natural approach to
consider, since both the index model [\citet{jiang.liu:2014}] and the varying
coefficient model [\citet{lian:2012}] cannot address this question; the former
is used to study interactions of components in a set of continuous
predictors, and the latter only allows one index variable. For model
selection and estimation, we apply a groupwise penalization method.
Moreover, most
existing high-dimensional nonparametric modeling papers [\citet%
{meier.van.Buhlmann:2009,Ravikumar.Liu.Lafferty.Wasserman:2009,
Huang.Horowitz.Wei:2010,lian:2012,wang.eta:2014}]
focus on variable selection and estimation. In this paper, after variable
selection, we also propose a simultaneous inferential tool to further test
the shape of the coefficient function for each selected variable, which
has not been studied in the previous works.

To this end, we aim to address questions (i)--(iii). Specifically,
for estimation and model selection, we apply a groupwise regularization
method based on a penalized quasi-likelihood criterion. The penalty is
imposed on the $L_{2}$ norm of the spline coefficients of the spline
estimators for $\alpha_{\ell}(\mathbf{\cdot})$. We establish the
asymptotic consistency of model selection and estimation for the proposed
group penalized estimators with the quasi-likelihood criterion in the
high-dimensional GACM (\ref{MOD:ADDITIVE}). We allow $p$ to grow with
$n$ at
an almost exponential order. Importantly, establishment of these
results is
technically more difficult than other work based on least squares,
since no
closed-form of the estimators exists from the penalized quasi-likelihood
method.

After selecting the important variables, the next question of interest is
what shapes the nonzero coefficient functions may have. Then we need to
provide an inferential tool to further check whether a coefficient function
has some specific parametric form. For example, when it is a constant
or a
linear function, the corresponding covariate has no or linear interaction
effects with another covariate, respectively. For global inference, we
construct simultaneous confidence bands (SCBs) for the nonparametric
additive functions based on a two-step estimation procedure. By using the
selected variables, we first propose a refined two-step spline
estimator for
the function of interest, which is proved to have a pointwise asymptotic
normal distribution and oracle efficiency. We then establish the bounds for
the SCBs based on the absolute maxima distribution of a Gaussian
process and
on the strong approximation lemma [\citet{Csorgo.Revesz:1981}]. Some other
related works on SCBs for nonparametric functions include \citet%
{hall:titterington:1988,hardle.marron:1991,claeskens.VanKeilegom:2003},
among others. We provide an asymptotic formula for the standard
deviation of
the spline estimator for the coefficient function, which involves unknown
population parameters to be estimated. The formula has somewhat complex
expressions and contains many parameters. Direct estimation therefore
may be
not accurate, particularly with the small or moderate sample sizes. As an
alternative, the bootstrap method provides us a reliable way to calculate
the standard deviation by avoiding estimating those population parameters.
We here apply the smoothed bootstrap method suggested by \citet{efron:2013},
which advocated that the method can improve coverage probability to
calculate the pointwise estimated standard deviations for the
estimators of
the coefficient functions. This method was originally proposed for
calculating the estimated standard deviation of the estimate of a parameter
of interest, such as the conditional mean. We extend this method to the case
of functional estimation. We demonstrate by simulation studies in
Section~\ref{sec7} that compared to the traditional resampling bootstrap method,
the smoothed bootstrap method can successfully improve the empirical
coverage rate.

The paper is organized as follows. Section~\ref{sec:method} introduces the
B-spline estimation procedure for the nonparametric functions,
describes the
adaptive group Lasso estimators and the initial Lasso estimators and
presents asymptotic results. Section~\ref{sec:twostep} describes the
two-step spline estimators and introduces the simultaneous confidence bands
and the bootstrap methods for calculating the estimated standard deviation.
Section~\ref{sec7} describes simulation studies, and Section~\ref{sec8}
illustrates the method through the analysis of an obesity data set from a
genome-wide association study. Proofs are in the \hyperref[app]{Appendix} and additional
supplementary material [\citet{ma.carroll.liang.xu2015}].

\section{Penalization based variable selection}\label{sec:method}

Let $(Y_{i},\mathbf{X}_{i}^{\mathrm{T}},\mathbf{T}_{i}^{\mathrm{T}})$, $
i=1,\ldots,n$, be random vectors that are independently and identically
distributed as $(Y,\mathbf{X}^{\mathrm{T}},\mathbf{T}^{\mathrm{T}})$, where
$%
\mathbf{X}_{i}=(X_{i1},\ldots,X_{id})^{\mathrm{T}}$ and
$\mathbf
{T}%
_{i}=(T_{i1},\ldots,T_{ip})^{\mathrm{T}}$. Write the negative
quasi-likelihood function $Q(\mu,y)=\int_{\mu}^{y}\{ (
y-\zeta ) /V(\zeta)\}\,d\zeta$. Estimation of the mean function
can be
achieved by minimizing the negative quasi-likelihood of the observed data
%
\begin{equation}
\sum_{i=1}^{n}Q\bigl\{g^{-1}\bigl
\{\eta(\mathbf{X}_{i},\mathbf{T}%
_{i})\bigr
\},Y_{i}\bigr\}. \label{EQ:QUASI}
\end{equation}

\subsection{Spline approximation}

We approximate the smooth functions $\alpha_{\ell k}(\cdot)$, $1\leq
k\leq
d$ and $1\leq\ell\leq p$ in (\ref{MOD:ADDITIVE}) by B-splines. As in most
work on nonparametric smoothing, estimation of the functions $\alpha
_{\ell
k}(\cdot)$ is conducted on compact sets. Without loss of generality, let
the compact set be $\mathcal{X}=[0,1]$. Let $G_{n}^{0}$ be the space of
polynomial splines of order $q\geq2$. We introduce a sequence of spline
knots
\[
t_{-q-1}=\cdots=t_{-1}=t_{0}=0<t_{1}<
\cdots<t_{N}<1=t_{N+1}=\cdots=t_{N+q},
\]
where $N\equiv N_{n}$ is the number of interior knots. In the
following, let
$J_{n}=N_{n}+q$. {For $0\leq j\leq N$, let $H_{j}=t_{j+1}-t_{j}$ be the
distance between neighboring knots and let $H=\max_{0\leq s\leq
N}H_{j}$. Following \citet{zhou.eta:1998},
to study asymptotic properties of the spline estimators for $\alpha
_{\ell
k}(\cdot)$, we assume that $\max_{0\leq j\leq N-1}|H_{j+1}-H_{j}|=o(N^{-1})$
and $H/\min_{0\leq j\leq N}H_{j}\leq M$, where $M>0$ is a predetermined
constant. Such an assumption is necessary for numerical implementation. In
practice, we can use the quantiles as the locations of the knots.}
Let $\{b_{j,k}(x_{k}):1\leq j\leq J_{n}\}^{\mathrm{T}}$ be the $q$th
order B spline basis functions given on page 87 of \citet{deboor:2001}.
For positive numbers $a_{n}$ and $b_{n}$, $a_{n}\asymp b_{n}$ means
that $%
\lim_{n\rightarrow\infty}a_{n}/b_{n}=c$, where $c$ is some nonzero finite
constant. For $1\leq j\leq J_{n}$, we adopt the centered B-spline functions
given in \citet{xue.yang:2006} such that $B_{j,k}(x_{k})=\sqrt{N}%
[b_{j,k}(x_{k})-\{E(b_{j,k})/E(b_{1,k})\}b_{1,k}(x_{k})]$, so that $%
E\{B_{j,k}(X_{k})\}=0$ and $\operatorname{var}\{B_{j,k}(X_{k})\}\asymp1$. Define the space
$G_{n}$ of additive spline functions as the linear space spanned by $B(%
\mathbf{x})=\{1,B_{j,k}(x_{k}),1\leq j\leq J_{n},1\leq k\leq d\}
^{\mathrm{T}%
} $, where $\mathbf{x}=(x_{1},\ldots,x_{d})^{\mathrm{T}}$. According
to the
result on page 149 of \citet{deboor:2001}, for $\alpha_{\ell k}(\cdot)$
satisfying condition (C3) in Appendix \ref{appa2} such that $\alpha_{\ell
k}^{(r-1)}(x_{k})\in C^{0,1}[0,1]$ for given integer $r\geq1$,
where $%
C^{0,1}[0,1]$ is the space of Lipschitz continuous functions on $[0,1]$
defined in Appendix \ref{appa2}, there is a function
%
\begin{equation}
\alpha_{\ell k}^{0}(x_{k})=\sum
_{j=1}^{J_{n}}\gamma_{j,\ell
k}B_{j}(x_{k})
\in G_{n}^{0}, \label{DEF:glzl}
\end{equation}
such that $\operatorname{sup}_{x_{k}\in{}[0,1]}|\alpha_{\ell
k}^{0}(x_{k})-\alpha_{\ell k}(x_{k})|=O(J_{n}^{-r})$. Then for every
$1\leq
\ell\leq p$, $\alpha_{\ell}(\mathbf{x})$ can be approximated well
by a
linear combination of spline functions in $G_{n}^{0}$, so that
%
\begin{equation}
\alpha_{\ell}(\mathbf{x})\approx\alpha_{\ell}^{0}(
\mathbf {x})=\gamma _{\ell0}+\sum_{k=1}^{d}
\sum_{j=1}^{J_{n}}\gamma _{j,\ell
k}B_{j,k}(x_{k})=B(
\mathbf{x})^{\mathrm{T}}{\bolds{\gamma}}_{\ell}, \label{EQ:spline}
\end{equation}
where ${\bolds{\gamma}}_{\ell}=(\gamma_{\ell0},%
{\bolds{\gamma}}_{\ell1}^{\mathrm{T}},\ldots,{\bolds{\gamma}}_{\ell
d}^{\mathrm{T}})^{\mathrm{T}}$, in which ${\bolds{\gamma}}_{\ell
k}=(\gamma
_{j,\ell k}:1\leq j\leq J_{n})^{\mathrm{T}}$. Thus the minimization
problem in
(\ref{EQ:QUASI}) is equivalent to finding $\widetilde{{\bolds{\gamma}}}^{0}=
(\widetilde{{\bolds{\gamma}}}_{\ell}^{0\mathrm{T}},1\leq\ell\leq p)^{\mathrm{T}}$ with $
\widetilde{{\bolds{\gamma}}}_{\ell}^{0}=(%
\widetilde{\gamma}_{\ell0}^{0},\widetilde{{\bolds{\gamma}}}_{\ell
1}^{0\mathrm{T}},\ldots,\widetilde{{\bolds{\gamma}}}_{\ell d}^{0
\mathrm{T}})^{\mathrm{T}}$ and $\widetilde{{\bolds{\gamma}}}_{\ell
k}^{0}=(\widetilde{\gamma}_{j,\ell k}^{0}:1\leq j\leq J_{n})^{\mathrm
{T}}$ to
minimize $\sum_{i=1}^{n}Q[g^{-1}\{\sum_{\ell=1}^{p}B(
\mathbf{X}_{i})^{\mathrm{T}}{\bolds{\gamma}}_{\ell
}T_{\ell
}\},Y_{i}]$. The components of the additive coefficients are estimated
by $%
\widetilde{\alpha}_{\ell k}^{0}(x_{k})=\sum_{j=1}^{J_{n}}\widetilde
{\gamma}%
_{j,\ell k}^{0}B_{j}(x_{k})=B(\mathbf{x})^{\mathrm{T}}\widetilde{%
{\bolds{\gamma}}}_{\ell k}^{0}$ and $\widetilde{\alpha
}_{\ell
0}^{0}=\widetilde{\gamma}_{\ell0}^{0}$.

\subsection{Adaptive group Lasso estimator}

We now describe the procedure for estimating and selecting the additive
coefficient functions by using the adaptive group Lasso. The estimators are
obtained by minimizing a penalized negative quasi-likelihood criterion. We
establish asymptotic selection consistency as well as the convergence rate
of the estimators to the true nonzero functions. For any vector
$\mathbf
{a}%
=(a_{1},\ldots,a_{s})^{\mathrm{T}}$, let its $L_{2}$ norm be
$\Vert
\mathbf{a}\Vert_{2}=\sqrt{a_{1}^{2}+\cdots+ a_{s}^{2}}$. For any
measurable $
L_{2}$-integrable function $\phi$ on $[0,1]^{d}$, define the $L_{2}$ norm
as $\Vert\phi\Vert^{2}=E\{\phi^{2}(\mathbf{X})\}$.

We are interested in identifying the significant components of the
vector $%
\mathbf{T}=(T_{1},\ldots,T_{p})^{\mathrm{T}}$. Let $s$, a fixed
number, be
the total number of nonzero $\alpha_{\ell}$'s and $I_{1}=\{\ell
:\Vert
\alpha_{\ell}\Vert\neq0,1\leq\ell\leq p\}$. {Let $I_{2}$ be the
complementary set of $I_{1}$; that is, $I_{2}=\{\ell:\alpha_{\ell
}(\cdot
)\equiv0,1\leq\ell\leq p\}$.} Recalling the approximation given in (%
\ref{EQ:spline}), ${\bolds{\gamma}}_{\ell}$ is
zero if
and only if each element of ${\bolds{\gamma}}_{\ell
}$ is
zero; that is, $\Vert{\bolds{\gamma}}_{\ell
}\Vert_{2}=0$.
We apply the adaptive group Lasso approach in \citet{Huang.Horowitz.Wei:2010}
for variable selection in model (\ref{MOD:ADDITIVE}). In order to identify
zero additive coefficients, we penalize the $L_{2}$ norm of the coefficients
${\bolds{\gamma}}_{\ell}$ for $1\leq\ell\leq p$.
Let $%
w_{n}=(w_{n1},\ldots,w_{np})^{\mathrm{T}}$ be a given vector of
weights, which needs to be chosen appropriately to achieve selection
consistency. Their choice will be discussed in Section~\ref{sec2.3}.
We consider the penalized negative quasi-likelihood
%
\begin{equation}
L_{n}({\bolds{\gamma}})=\sum_{i=1}^{n}Q
\Biggl[g^{-1}\Biggl\{%
\sum_{\ell=1}^{p}B^{\mathrm{T}}(
\mathbf{X}_{i}){\bolds{\gamma}}_{\ell}T_{\ell}\Biggr
\},Y_{i}\Biggr]+n\lambda _{n}\sum
_{\ell=1}^{p}w_{n\ell}\Vert{\bolds{
\gamma}}_{\ell}\Vert_{2}, \label{EQ:LASSO}
\end{equation}
where $\lambda_{n}$ is a regularization parameter controlling the
amount of
shrinkage. The estimator ${\widehat{{\bolds{\gamma}}}}%
=({\widehat{{\bolds{\gamma}}}}_{1}^{\mathrm
{T}%
},\ldots,{\widehat{{\bolds{\gamma}}}}_{p}^{\mathrm{T}})^{
\mathrm{T}}$ is obtained by minimizing (\ref{EQ:LASSO}). {Minimization of
(\ref%
{EQ:LASSO}) is solved by local quadratic approximation as adopted by
\citet%
{fan.li:2001}.}

For $\ell=1,\ldots,p$, the $\ell$th additive coefficient function is
estimated by
\[
\widehat{\alpha}_{\ell}(\mathbf{x})=\widehat{\gamma}_{\ell0}+%
\sum_{k=1}^{d}\sum
_{j=1}^{J_{n}}\widehat{\gamma}_{j,\ell
k}B_{j,k}(x_{k})=B^{\mathrm{T}}(
\mathbf{x}){\widehat{%
{\bolds{\gamma}} }}_{\ell}.
\]

We will make the following two assumptions on the order requirements of the
tuning parameters. Write $w_{n,I_{1}}=(w_{n \ell}:\ell\in I_{1})$.

\begin{assumption}
\label{Assumption1} $J_{n}^{2}\{n\log(n)\}^{-1}\rightarrow0$ and
$\lambda
_{n}\Vert w_{n,I_{1}}\Vert_{2}\rightarrow0$, as $n\rightarrow\infty$.
\end{assumption}

\begin{assumption}
\label{Assumption2} $n\lambda_{n}\Vert w_{n,I_{1}}\Vert
_{2}+n^{1/2}J_{n}^{1/2}\sqrt{\log( pJ_{n}) }+nJ_{n}^{-r}=o( n\lambda
_{n}w_{n \ell}) $, for all $\ell\in I_{2}$.
\end{assumption}

The following theorem presents the selection consistency and estimation
properties of
the adaptive group Lasso estimators.

\begin{theorem}
\label{THM:weakconsistencyLASSORE} Under conditions \textup{(C1)--(C5)} in the
\hyperref[app]{Appendix} and Assumptions~\ref{Assumption1} and
\ref{Assumption2}: \textup{(i)}
as $%
n\rightarrow\infty$, $P(\llVert \widehat{\alpha}_{\ell}\rrVert
>0,\ell\in I_{1}\mbox{ and }\llVert \widehat{\alpha}_{\ell
}\rrVert
=0,\ell\in I_{2})\rightarrow1$, and \textup{(ii)} $\Vert\widehat{\alpha
}_{\ell
}-\alpha_{\ell}\Vert=O_{p}(\lambda_{n}\Vert w_{n,I_{1}}\Vert
_{2}+n^{-1/2}J_{n}^{1/2}+J_{n}^{-r}),\ell\in I_{1}$.
\end{theorem}

\subsection{Choice of the weights}\label{sec2.3}

We now discuss how to choose the weights used in (\ref{EQ:LASSO})
based on
the initial estimates. For low-dimensional data settings with $p<n$, an
unpenalized estimator such as least squares estimator [\citet{zou:2006}] can
be used as an initial estimate. For high-dimensional settings with
$p\gg n$,
it has been discussed [\citet{meier.Buhlmann:2007}] that the Lasso estimator
is a more appropriate choice. Following \citet
{Huang.Horowitz.Wei:2010}, we
obtain an initial estimate with the group Lasso by minimizing
\[
L_{n1}({\bolds{\gamma}})=\sum_{i=1}^{n}Q
\Biggl[g^{-1}\Biggl\{%
\sum_{\ell=1}^{p}B(
\mathbf{X}_{i})^{\mathrm{T}}{\bolds{\gamma}}_{\ell}T_{\ell}
\Biggr\},Y_{i}\Biggr]+n\lambda _{n1}\sum
_{\ell=1}^{p}\Vert{\bolds{\gamma}} 
_{\ell}
\Vert_{2},
\]
with respect to ${\bolds{\gamma}}=({\bolds{\gamma}}_{1}^{\mathrm{T}},\ldots,
{\bolds{\gamma}}_{p}^{\mathrm{T}})^{\mathrm{T}}$. Denote the resulting estimators
by $%
\widetilde{{\bolds{\gamma}}}=(\widetilde{{\bolds{\gamma}}}_{1}^{\mathrm{T}},\ldots,\widetilde{{\bolds{\gamma}}}_{p}^{\mathrm
{T}})^{\mathrm{T}}$. Let $\widetilde{I}%
_{1}=\{\ell:\Vert\widetilde{{\bolds{\gamma}}}_{\ell
}\Vert_{2}\neq0,1\leq\ell\leq p\}$, and let $\widetilde{s}$ be the
number of elements in $\widetilde{I}_{1}$.

Under conditions (C1)--(C5) in the \hyperref[app]{Appendix}, and when $\lambda
_{n1}\geq
Cn^{-1/2}\times\break J_{n}^{1/2}  \sqrt{\log(pJ_{n})}$ for a sufficiently large
constant $%
C$, we have: (i) the number of estimated nonzero functions are bounded;
that is,
as $n\rightarrow\infty$, there exists a constant $1<C_{1}<\infty$ such
that $P(\widetilde{s}\leq C_{1}s)\rightarrow1$; (ii) if $\lambda
_{n1}\rightarrow0$,
then $P(\Vert\widetilde{{\bolds{\gamma}}}_{\ell}\Vert_{2}>0$ for all
\mbox{$l\in I_{1})\rightarrow1$}; (iii) $\Vert\widetilde{{\bolds{\gamma}}}-{\bolds{\gamma}}\Vert
_{2}=O_{p}(\lambda
_{n1}+n^{-1/2}J_{n}^{1/2}+J_{n}^{-r})$. We refer to Theorems~1(i) and (ii)
of \citet{Huang.Horowitz.Wei:2010} for the proofs of (i) and (ii), and
Theorem~\ref{THM:weakconsistencyLASSORE} in our paper for the proof of (iii).

The weights we use are $w_{n\ell}=\Vert\widetilde{{\bolds{\gamma}}}_{\ell}\Vert_{2}^{-1}$, if
$\Vert\widetilde{{\bolds{\gamma}}}_{\ell}\Vert_{2}>0$; $w_{n\ell}=\infty$,
if $%
\Vert\widetilde{{\bolds{\gamma}}}_{\ell}\Vert_{2}=0$.

\begin{remark}
\label{Remark1} Assumptions \ref{Assumption1} and \ref{Assumption2}
give the
order requirements of $J_{n}$ and~$\lambda_{n}$. Based on the condition
that $J_{n}^{2}\{n\log(n)\}^{-1}\rightarrow0$ given in Assumption~\ref
{Assumption1}, we need $J_{n}\ll\{n\log(n)\}^{1/2}$, where $a_{n}\ll b_{n}$
denotes that $a_{n}/b_{n}=o ( 1 ) $ for any positive numbers $a_{n}$
and $b_{n}$, and $\lambda_{n}$ needs to satisfy
$n^{-1/2}J_{n}^{1/2}\times\break \sqrt{%
\log(pJ_{n})} \{\min_{\ell\in I_{2}}(w_{n\ell})\}
^{-1}\ll\lambda
_{n}\ll1$. From the above theoretical properties of the group Lasso
estimators, we know that, with probability approaching 1, $\Vert
\widetilde{{\bolds{\gamma}}}_{\ell}\Vert_{2}>0$ for nonzero
components, and then the corresponding weights $w_{n\ell}$ are bounded away
from $0$ and infinity for $\ell\in I_{1}$. By defining $0\cdot\infty=0$,
the components not selected by the group Lasso are not included in the
adaptive group Lasso procedure. Let $J_{n}\asymp n^{1/(2r+1)}$, so that
$%
J_{n}$ has the optimal order for spline regression. If $p=\exp
[o\{n^{2r/(2r+1)}\}]$, then $n^{-1/2}J_{n}^{1/2}\sqrt{\log(pJ_{n})}%
\rightarrow0$. This means the dimension $p$ can diverge with the sample
size at an almost exponential rate.
\end{remark}

\subsection{Selection of tuning parameters}
\label{sec:computation}

Tuning parameter selection always plays an important role in model and
variable selection. An underfitted model can lead to severely biased
estimation, and an overfitted model can seriously degrade the estimation
efficiency. Among different data-driven methods, the Bayesian information
criterion (BIC) tuning parameter selector has been shown to be able to
identify the true model consistently in the fixed dimensional setting %
[\citet{wang.li.tsai:2007}]. In the high-dimensional setting, an extend BIC
(EBIC) and a generalized information criterion have been proposed by
\citet%
{chen.chen:2008} and \citet{fan.tang:2013}, respectively. In this
paper, we
adopt the EBIC method [\citet{chen.chen:2008}] to select the tuning
parameter $%
\lambda_{n}$ in (\ref{EQ:LASSO}). Specifically, the EBIC$(\lambda
_{n})$ is defined as
\[
2\sum_{i=1}^{n} \Biggl( Q\Biggl[
g^{-1}\Biggl\{ \sum_{\ell=1}^{p}B(
\mathbf{X}_{i}) ^{\mathrm{T}}\widehat{{\bolds{\gamma}
}}_{\ell}T_{i\ell}\Biggr\},Y_{i}\Biggr] \Biggr) +
s^{\ast}( 1+dJ_{n}) \operatorname{log}( n) +2\nu\operatorname{log}\pmatrix{p
\cr
s^{\ast}},
\]
where $( \widehat{{\bolds{\gamma} }}_{\ell})
_{\ell
=1}^{p}$ is the minimizer of (\ref{EQ:LASSO}) for a given $\lambda
_{n}$, $%
s^{\ast}$ is the number of nonzero estimated functions $(\widehat
{\alpha}%
_{\ell}) _{\ell=1}^{p}$ and $0\leq\nu\leq1$ is a constant. Here we use
$\nu=0.5$. When $\nu=0$, the EBIC is ordinary BIC.

We use cubic B-splines for the nonparametric function estimation, so
that $%
q=4$. In the penalized estimation procedure, we let the number of interior
knots $N=\lfloor cn^{1/( 2q+1) }\rfloor$ satisfy the optimal order,
where $%
\lfloor a\rfloor$ denotes the largest integer no greater than $a$ and $c$
is a constant. In the simulations, we take $c=2$.

\section{Inference and the bootstrap smoothing procedure}\label{sec:twostep}

\subsection{Background}\label{sec3.11}

After model selection, our next step is to conduct statistical
inference for
the coefficient functions of those important variables. We will
establish a
simultaneous confidence band (SCB) based on a two-step estimator for global
inference. An asymptotic formula of the SCB will be provided based on the
distribution of the maximum value of the normalized deviation of the spline
functional estimate. To improve accuracy, we calculate the estimated
standard deviation in the SCB by using the nonparametric bootstrap smoothing
method as discussed in \citet{efron:2013}. For specificity, we focus on the
construction of $\alpha_{\ell1}(x_{1})$, with $\alpha_{\ell k}(x_{k})$
for $k\geq2$ defined similarly, for $\ell\in\widehat{I}_{1}$, where $
\widehat{I}_{1}=\{\ell:\Vert\widehat{\alpha}_{\ell}\Vert\neq
0,1\leq
\ell\leq p\}$.

Although the one-step penalized estimation in Section~\ref{sec:method} can
quickly identify nonzero coefficient functions, no asymptotic distribution
is available for the resulting estimators. Thus we construct the SCB based
on a refined two-step spline estimator for $\alpha_{\ell1}(x_{1})$, which
will be shown to have the oracle property that the estimator of $\alpha
_{\ell1}(x_{1})$ has the same asymptotic distribution as the univariate
oracle estimator obtained by pretending 
that $\alpha_{\ell0}$ and $\alpha_{\ell k}(X_{k})$ for $\ell\in
\widehat{%
I}_{1}$, $k\geq2$ and $\alpha_{\ell}(\mathbf{X})$ for $\ell\notin
\widehat{I}_{1}$ are known. See \citet{horowitz.mammen:2004,Horowitz.Klemela.Mammen:2006,liu.yang.hardle:2013} for kernel-based two-step
estimators in generalized additive models, which also have the oracle
property but are not as computationally efficient as the two-step spline
method. We next introduce the oracle estimator and the proposed two-step
estimator before we present the SCB.

\subsection{Oracle estimator}\label{sec3.1b}

In the following, we describe the oracle estimator of $\alpha_{\ell
1}(x_{1})$. We rewrite model (\ref{MOD:ADDITIVE}) as
%
\begin{eqnarray}
\label{eq0A} \mu(\mathbf{X},\mathbf{T}) &=&g^{-1}\bigl\{\eta(
\mathbf{X},\mathbf{T})\bigr\}
\nonumber
\\[-8pt]
\\[-8pt]
\nonumber
&=&\sum_{\ell
\in\widehat{I}_{1}}
\alpha_{\ell
1}(X_{1})T_{\ell}+\sum
_{\ell\in\widehat{I}_{1}}\biggl\{\alpha _{\ell
0}
+\sum_{k\geq2}\alpha_{\ell k}(X_{k})
\biggr\}T_{\ell
}+\sum_{\ell\notin\widehat{I}_{1}}
\alpha_{\ell}(\mathbf{X}%
)T_{\ell}.
\end{eqnarray}
By assuming that $\alpha_{\ell0}$ and $\alpha_{\ell k}(X_{k})$ for
$\ell
\in\widehat{I}_{1}$, $k\geq2$ and $\alpha_{\ell}(\mathbf{X})$ for
$\ell
\notin\widehat{I}_{1}$ are known, estimation in (\ref{eq0A})
involves only
the nonparametric functions $\alpha_{\ell1}(X_{1})$ of a scalar
covariate $%
X_{1}$. It will be shown in Theorem~\ref{THM:asymptotics} that the estimator
achieves the univariate optimal convergence rate when the optimal order for
the number of knots is applied. We estimate $\alpha_{1}(x_{1})=\{
\alpha
_{\ell1}(x_{1}),\ell\in\widehat{I}_{1}\}^{\mathrm{T}}$ by minimizing the
negative quasi-likelihood function as follows. Denote the oracle estimator
by $\widehat{\alpha}_{\ell1}^{\mathrm{OR}}(x_{1})=B_{1}^{\mathcal{S}%
}(x_{1})^{\mathrm{T}}\widehat{{\bolds{\gamma}}}_{\ell
1}^{\mathrm
{OR}}$%
, where $\widehat{{\bolds{\gamma}}}_{\ell1}^{\mathrm{OR}}$ is
defined directly below, $B_{1}^{\mathcal{S}}(x_{1})=\{
B_{j,1}^{\mathcal
{S}%
}(x_{1}),1\leq j\leq J_{n}^{\mathcal{S}}\}$ where $B_{j,1}^{\mathcal
{S}%
}(x_{1})$ is the centered B-spline function defined in the same way as $
B_{j,1}(x_{1})$ in Section~\ref{sec:method}, but with $N^{\mathcal{S}%
}=N_{n}^{\mathcal{S}}$ interior knots and $J_{n}^{\mathcal{S}}=N_{n}^{
\mathcal{S}}+q$. Rates of increase for $J_{n}^{\mathcal{S}}$ are described
in Assumptions \ref{Assumption3} and \ref{Assumption4} below. Let
$\alpha
_{\ell,-1}(\mathbf{X}_{i})=\alpha_{\ell0}+\sum_{k\geq
2}\alpha
_{\ell k}(X_{ik})$. Then ${\widehat{{\bolds{\gamma}}}}_{,1}^{%
\mathrm{OR}}=\{(\widehat{{\bolds{\gamma}}}_{\ell
1}^{\mathrm{OR%
}})^{\mathrm{T}},\ell\in\widehat{I}_{1}\}^{\mathrm{T}}$ is obtained by
minimizing the negative quasi-likelihood
%
\begin{eqnarray}
\label{EQ:LnOR} L_{n}^{\mathrm{OR}}({\bolds{\gamma}}_{,1})
&=&\sum_{i=1}^{n}Q\biggl[g^{-1}
\biggl\{\sum_{\ell\in\widehat{I}%
_{1}}B_{1}^{\mathcal{S}}(X_{i1})^{\mathrm{T}}{
\bolds{\gamma}}_{\ell
1}T_{i\ell}
\nonumber
\\[-8pt]
\\[-8pt]
\nonumber
&&{}+\sum_{\ell\in\widehat{I}_{1}}\alpha_{\ell,-1}(\mathbf{X}
_{i})T_{i\ell}+\sum_{\ell\notin\widehat{I}_{1}}
\alpha _{\ell}(%
\mathbf{X}_{i})T_{i\ell}
\biggr\},Y_{i}\biggr],
\end{eqnarray}
where ${\bolds{\gamma}}_{,1}=\{(\bolds{\gamma}_{\ell
1})^{\mathrm{T}},\ell\in\widehat{I}_{1}\}^{\mathrm{T}}$. Similarly,
the oracle
estimator of $\bolds{\alpha}_{0}=\{\alpha_{\ell0},\ell
\in
\widehat{I}_{1}\}^{\mathrm{T}}$, which is denoted as
$\widehat{\bolds{\alpha}}_{0}^{\mathrm{OR}}=\{\widehat{\alpha}_{\ell0}^{\mathrm
{OR}},\ell\in
\widehat{I}_{1}\}^{\mathrm{T}}=\{\widehat{\gamma}_{\ell0}^{\mathrm
{OR}},\ell
\in\widehat{I}_{1}\}^{\mathrm{T}}$, is obtained by minimizing
$
L_{n}^{\mathrm{OR}}({\bolds{\gamma}}_{,0})=\sum%
_{i=1}^{n}Q[g^{-1}\{\sum_{\ell\in\widehat
{I}_{1}}\gamma
_{\ell0}T_{i\ell}+\break \sum_{\ell\in\widehat{I}_{1}}\alpha
_{\ell
,-0}(\mathbf{X}_{i})T_{i\ell}+\sum_{\ell\notin\widehat{I}%
_{1}}\alpha_{\ell}(\mathbf{X}_{i})T_{i\ell}\},Y_{i}]$,
where ${\bolds{\gamma}}_{,0}=({\bolds{\gamma}}_{\ell
0},\ell\in\widehat{I}_{1})$ and\break $\alpha_{\ell,-0}(\mathbf{X}%
_{i})=\sum_{k=1}^{d}\alpha_{\ell k}(X_{ik})$.

\subsection{Initial estimator}\label{sec3.2}

The oracle estimator is infeasible because it assumes knowledge of the other
functions. In order to obtain the two-step estimators of $\alpha_{\ell
1}(x_{1})$ for $\ell\in\widehat{I}_{1}$, we first need initial estimators
for $\alpha_{\ell0}$ and $\alpha_{\ell k}(x_{k})$ for $k\geq2$ and
$\ell
\in\widehat{I}_{1}$, denoted as $\widehat{\alpha}_{\ell0}^{\mathrm
{ini}}=%
\widehat{\gamma}_{\ell0}^{\mathrm{ini}}$ and $\widehat{\alpha
}_{\ell
k}^{%
\mathrm{ini}}(x_{k})=B_{k}^{\mathrm{ini}}(x_{k})^{\mathrm{T}}\widehat{{\bolds{\gamma}}}_{\ell k}^{\mathrm{ini}}$, where
$B_{k}^{\mathrm
{ini}%
}(x_{k})=\{B_{j,k}^{\mathrm{ini}}(x_{k}):1\leq j\leq J_{n}^{\mathrm
{ini}}\}
^{%
\mathrm{T}}$ and $B_{j,k}^{\mathrm{ini}}(x_{k})$ are B-spline functions with
the number of interior knots $N_{n}^{\mathrm{ini}}$ and $J_{n}^{\mathrm{ini}
}=N_{n}^{\mathrm{ini}}+q$. Rates of increase for $J_{n}^{\mathrm{ini}}$ are
described in Assumptions \ref{Assumption3} and \ref{Assumption4}
below. We
need an undersmoothed procedure in the first step, so that the approximation
bias can be reduced, and the difference between the two-step and oracle
estimators is asymptotically negligible. We obtain ${\widehat{{\bolds{\gamma}}}}_{\widehat{I}_{1}}^{\mathrm{ini}}=\{
({\widehat{{\bolds{\gamma}}}}_{\ell}^{\mathrm{ini}})^{\mathrm
{T}}:\ell
\in\widehat{I}_{1}\}^{\mathrm{T}}$, where ${\widehat{{\bolds{\gamma}}}}_{\ell}^{\mathrm{ini}}=
\{\widehat{{\bolds{\gamma}}}_{\ell0}^{\mathrm{ini}},(\widehat{{\bolds{\gamma}}}_{\ell
k}^{\mathrm{ini}})^{\mathrm{T}}\}^{\mathrm{T}}$, by minimizing the negative
quasi-likelihood $\sum_{i=1}^{n}Q[g^{-1}\{\sum_{\ell
\in
\widehat{I}_{1}}B(\mathbf{X}_{i})^{\mathrm{T}}{\bolds{\gamma}}_{\ell}T_{\ell}\},Y_{i}]$. The adaptive group Lasso penalized
estimator ${\widehat{{\bolds{\gamma}}}}_{\widehat
{I}_{1}}=\{(%
{\widehat{{\bolds{\gamma}}}}_{\ell})^{\mathrm
{T}}:\ell
\in
\widehat{I}_{1}\}^{\mathrm{T}}$ obtained in Section~\ref{sec:method}
can also
be used as the initial estimator. We, however, refit the model with the
selected variables and obtain the initial estimator ${\widehat{{\bolds{\gamma}}}}_{\widehat{I}_{1}}^{\mathrm{ini}}$ in
order to
improve estimation accuracy in high-dimensional data settings.

\subsection{Final estimator}\label{sec3.3}

In the second step, we construct the two-step estimator of $\alpha
_{\ell1}$
for $\ell\in\widehat{I}_{1}$. We replace $\alpha_{\ell0}$ and
$\alpha
_{\ell k}(X_{k})$ by the initial estimators $\widehat{\alpha}_{\ell
0}^{%
\mathrm{ini}}$ and $\widehat{\alpha}_{\ell k}^{\mathrm{ini}}(X_{k})$ for
$\ell
\in\widehat{I}_{1}$ and $k\geq2$ and replace $\alpha_{\ell
}(\mathbf{X})$
for $\ell\notin\widehat{I}_{1}$ by $\widehat{\alpha}_{\ell
}(\mathbf
{X}%
)=0 $. Let $\widehat{\alpha}_{\ell,-1}^{\mathrm{ini}}(\mathbf{X}_{i})=
\widehat{\alpha}_{\ell0}^{\mathrm{ini}}+\sum_{k\geq2}\widehat{
\alpha}_{\ell k}^{\mathrm{ini}}(X_{ik})$. Denote the two-step spline
estimator of $\alpha_{\ell1}(x_{1})$ as $\widehat{\alpha}_{\ell1}^{
\mathcal{S}}(x_{1})=B_{1}^{\mathcal{S}}(x_{1})^{\mathrm{T}}\widehat{{\bolds{\gamma}}}_{\ell1}^{\mathcal{S}}$ with
${\widehat{{\bolds{\gamma}}}}_{,1}^{\mathcal{S}}=\{(\widehat{{\bolds{\gamma}}}_{\ell
1}^{\mathcal{S}})^{\mathrm{T}},\ell\in\widehat{I}_{1}\}^{\mathrm{T}}$
minimizing
%
\begin{eqnarray}
\label{EQ:Lns}  L_{n}^{\mathcal{S}}({\bolds{\gamma}}_{,1})&=&
\sum%
_{i=1}^{n}Q\biggl[g^{-1}\biggl\{\sum
_{\ell\in\widehat
{I}_{1}}B_{1}^{%
\mathcal{S}}(X_{i1})^{\mathrm{T}}{
\bolds{\gamma}}_{\ell1}T_{i\ell}
\nonumber
\\[-8pt]
\\[-8pt]
\nonumber
&&{} +\sum_{\ell\in\widehat{I}_{1}}\widehat{\alpha }%
_{\ell,-1}^{\mathrm{ini}}(
\mathbf{X}_{i})T_{i\ell}+\sum_{\ell
\notin\widehat{I}_{1}}
\widehat{\alpha}_{\ell}(\mathbf {X}_{i})T_{i\ell
}
\biggr\},Y_{i}\biggr].
\end{eqnarray}
Then the two-step of $\alpha_{\ell0}$, denoted as $\widehat{\alpha
}_{\ell
0}^{\mathcal{S}}=\widehat{\gamma}_{\ell0}^{\mathcal{S}}$, is
obtained in
the same way as $\widehat{\alpha}_{\ell0}^{\mathrm{OR}}$ by replacing
$%
\alpha_{\ell,0}(\mathbf{X}_{i})$ with $\widehat{\alpha}_{\ell
,0}^{\mathrm{%
ini}}(\mathbf{X}_{i})=\sum_{k=1}^{d}\widehat{\alpha}_{\ell
k}^{%
\mathrm{ini}}(X_{ik})$ for $\ell\in\widehat{I}_{1}$ and replacing
$\alpha
_{\ell}(\mathbf{X}_{i})$ with $\widehat{\alpha}_{\ell}(\mathbf{X}_{i})=0$
for $\ell\notin\widehat{I}_{1}$. Let $\widehat{\bolds{\alpha}}_{0}^{\mathcal{S}}=\{\widehat{\alpha}_{\ell0}^{\mathcal
{S}},\ell
\in\widehat{I}_{1}\}^{\mathrm{T}}$.

\subsection{Asymptotic normality and uniform oracle efficiency}\label{sec3.4}

We now establish the asymptotic normality and uniform oracle efficiency for
the oracle and final estimators. Let $Z_{ij\ell,1}=B_{j,1}^{\mathcal{S}}(
X_{i1}) T_{i\ell}$ and $Z_{i,1}=( Z_{ij\ell,1},1\leq j\leq
J_{n}^{\mathcal{S}%
},\ell\in\widehat{I}_{1}) ^{\mathrm{T}}$. Let $s^{\ast}$ be the
number of
elements in $\widehat{I}_{1}$. By Theorem~\ref
{THM:weakconsistencyLASSORE}, $%
P( s^{\ast}=s) \rightarrow1$. For simplicity of notation, denote
$\sigma
_{i}^{2}=\sigma^{2}( \mathbf{X}_{i},\mathbf{T}_{i}) $ and $\eta
_{i}=\eta(
\mathbf{X}_{i},\mathbf{T}_{i}) $. Define ${s^{\ast}\times s^{\ast
}J_{n}^{\mathcal{S}}}$
matrix $B^{\mathcal{S}}( x_{1})$ as
\[
\left[ %
\matrix{ B_{1,1}^{\mathcal{S}}(
x_{1}) & \cdots& B_{J_{n}^{\mathcal
{S}},1}^{\mathcal{%
S}}( x_{1})
& 0 & \cdots& 0 & 0 & \cdots& 0
\cr
\vdots& \vdots& \vdots& \vdots& \vdots& \vdots& \vdots& \vdots& \vdots
\cr
0 & \cdots& 0 & 0 & \cdots& 0 & B_{1,1}^{\mathcal{S}}(
x_{1}) & \cdots& B_{J_{n}^{\mathcal{S}},1}^{\mathcal{S}}(
x_{1})}
 \right].
\]

To establish the asymptotic distribution of the two-step estimator, in
addition to Assumptions 1 and 2 given in Section~\ref{sec:method}, we make
the following two assumptions on the number of basis functions $J_{n}^{%
\mathcal{S}}$ and $J_{n}^{\mathrm{ini}}$:

\begin{assumption}
\label{Assumption3} (i) $s^{\ast}(J_{n}^{\mathcal{S}})^{2}\{n\log
(n)\}^{-1}=o(1)$ and $s^{\ast}(J_{n}^{\mathcal{S}})^{-r}=o(1)$, and
(ii) $%
n(\log n)^{-1}(J_{n}^{\mathcal{S}}J_{n}^{\mathrm{ini}})^{-1}\rightarrow
\infty
$, as $n\rightarrow\infty$.
\end{assumption}

\begin{assumption}
\label{Assumption4} $( n/J_{n}^{\mathcal{S}}) ^{1/2}( J_{n}^{\mathrm{ini}})
^{-r}\rightarrow0$, as $n\rightarrow\infty$.
\end{assumption}

First we describe the asymptotic normality of the oracle estimator
$\widehat{%
\alpha}_{\ell1}^{\mathrm{OR}}(x_{1})$ of $\alpha_{\ell1}(x_{1})$. Let
$%
\widehat{\alpha}_{1}^{\mathrm{OR}}(x_{1})=\{\widehat{\alpha}_{\ell
1}^{\mathrm{%
OR}}(x_{1}),\ell\in\widehat{I}_{1}\}^{\mathrm{T}}$. Let
$b_{1}(x_{1})=E\{%
\widehat{\alpha}_{1}^{\mathrm{OR}}(x_{1})|\mathbb{X},\mathbb{T}\}$
and $%
b_{\ell1}(x_{1})=E\{\widehat{\alpha}_{\ell1}^{\mathrm
{OR}}(x_{1})|\mathbb{X}%
,\mathbb{T}\}$, for $\ell\in\widehat{I}_{1}$, where $(\mathbb
{X},\mathbb{T}%
)=(\mathbf{X}_{i},\mathbf{T}_{i})_{i=1}^{n}$.

\begin{theorem}\label{THM:asymptotics}
Under conditions \textup{(C1)--(C5)} and Assumption~\ref%
{Assumption3}\textup{(i)}, for any vector $\mathbf{a}\in R^{s^{\ast}}$ with
$\Vert
\mathbf{a}\Vert_{2}=1$, for any $x_{1}\in{}[0,1]$,
$\mathbf{a}^{\mathrm{T}}\sigma_{n}^{-1}(x_{1})\{\widehat{\alpha
}_{1}^{\mathrm{%
OR}}(x_{1})-b_{1}(x_{1})\}\rightarrow N(0,1)$, where
%
\begin{equation}
\sigma_{n}^{2}(x_{1})=B^{\mathcal{S}}(x_{1})
\Biggl[\sum%
_{i=1}^{n}Z_{i,1}Z_{i,1}^{\mathrm{T}}
\bigl\{\dot{g}^{-1}(\eta _{i})\bigr\}^{2}\Big/
\sigma_{i}^{2}\Biggr]^{-1}B^{\mathcal{S}}(x_{1})^{\mathrm{T}},
\label{DEF:sig2nx1}
\end{equation}
where $\dot{g}^{-1}(\eta_{i})$ is the first-order derivative of $%
g^{-1}(\eta_{i})$ with respect to $\eta_{i}$, and
\[
\sum_{\ell\in\widehat{I}_{1}}\bigl\Vert\widehat{\alpha}_{\ell
1}^{%
\mathrm{OR}}-b_{\ell1}
\bigr\Vert^{2}=O_{p}\bigl(s^{\ast}J_{n}^{\mathcal{S}%
}n^{-1}
\bigr),\qquad\sum_{\ell\in I_{1}}\Vert b_{\ell1}-
\alpha_{\ell
1}\Vert^{2}=O_{p}\bigl\{
\bigl(s^{\ast}\bigr)^{2}\bigl(J_{n}^{\mathcal{S}}
\bigr)^{-2r}\bigr\}.
\]
Thus for $\ell\in\widehat{I}_{1}$, $\sigma_{n1}^{-1}(x_{1})\{
\widehat
{%
\alpha}_{\ell1}^{\mathrm{OR}}(x_{1})-b_{\ell1}(x_{1})\}\rightarrow N(0,1)$, where
%
\begin{equation}
\sigma_{n1}^{2}(x_{1})=\mathbf{e}_{\ell}^{\mathrm{T}}
\sigma _{n}^{2}(x_{1})%
\mathbf{e}_{\ell}, \label{EQ:sign1}
\end{equation}
and $\mathbf{e}_{\ell}$ is the $s^{\ast}$-dimensional vector with the
$\ell$th element $1$ and other elements $0$, and $\Vert\widehat{\bolds{\alpha}}_{0}^{\mathrm{OR}}-
\bolds{\alpha}_{0}\Vert_{2}=O_{p}(%
\sqrt{s^{\ast}/n})$.
\end{theorem}

The next result shows the uniform oracle efficiency of the two-step
estimator that the difference between the two-step estimator $\widehat{
\alpha}_{\ell1}^{\mathcal{S}}(x_{1})$ and oracle estimator $\widehat
{%
\alpha}_{\ell1}^{\mathrm{OR}}(x_{1})$ is uniformly asymptotically
negligible, and thus the two-step estimator is oracle in the sense that it
has the same asymptotic distribution as the oracle estimator. Let
$\widehat{%
\alpha}_{1}^{\mathcal{S}}(x_{1})=\{\widehat{\alpha}_{\ell
1}^{\mathcal
{S}%
}(x_{1}),\ell\in\widehat{I}_{1}\}^{\mathrm{T}}$.

\begin{theorem}\label{THM:oracle}
Under conditions \textup{(C1)--(C5)} in the \hyperref[app]{Appendix} and
Assumptions~\ref{Assumption1}--\ref{Assumption3},
\[
\sup_{x_{1}\in{}[0,1]}\bigl\Vert\widehat{\alpha }_{1}^{\mathcal{S}%
}(x_{1})-
\widehat{\alpha}_{1}^{\mathrm{OR}}(x_{1})
\bigr\Vert_{\infty
}=O_{p}\bigl\{\bigl(n^{-1}\log n
\bigr)^{1/2}+\bigl(J_{n}^{\mathrm{ini}}\bigr)^{-r}
\bigr\},
\]
$\Vert\widehat{\bolds{\alpha}}_{0}^{\mathcal
{S}}-\widehat{\bolds{\alpha}}_{0}^{\mathrm{OR}}\Vert
_{2}=o_{p}(n^{-1/2})$, and
furthermore under Assumption~\ref{Assumption4},
\[
\sup_{x_{1}\in{}[0,1]}\bigl|\mathbf{a}^{\mathrm{T}}\sigma
_{n}^{-1}(x_{1})\bigl\{\widehat{
\alpha}_{1}^{\mathcal{S}}(x_{1})-\widehat {\alpha
}_{1}^{\mathrm{OR}}(x_{1})\bigr\}\bigr|=o_{p}(1),
\]
for any vector $\mathbf{a}\in R^{s^{\ast}}$ with $\Vert\mathbf
{a}\Vert
_{2}=1$ and $\sigma_{n}^{2}(x_{1})$ given in (\ref{DEF:sig2nx1}). Hence,
for any $x_{1}\in{}[0,1]$, $\mathbf{a}^{\mathrm{T}}\sigma
_{n}^{-1}(x_{1})\{\widehat{\alpha}_{1}^{\mathcal
{S}}(x_{1})-b_{1}(x_{1})\}%
\rightarrow N(0,1)$.
\end{theorem}

\begin{remark}
\label{Remark2} Under Assumptions \ref{Assumption1} and \ref{Assumption2},
by Theorem~\ref{THM:weakconsistencyLASSORE}, with probability
approaching $1$%
, $s^{\ast}=s$, which is a fixed number. In the second step, by
letting $%
J_{n}^{\mathcal{S}}\asymp n^{1/(2r+1)}$, the nonparametric functions
$\alpha
_{\ell1}$ for $\ell\in\widehat{I}_{1}$ are approximated by spline
functions with the optimal number of knots. By the conditions that $%
(n/J_{n}^{\mathcal{S}})(J_{n}^{\mathrm{ini}})^{-1}\rightarrow0$ and
$n(\log n)^{-1}(J_{n}^{\mathcal{S}}J_{n}^{\mathrm
{ini}})^{-1}\rightarrow
\infty$ given in Assumptions \ref{Assumption3} and \ref
{Assumption4}, $
J_{n}^{\mathrm{ini}}$ needs to satisfy $n^{1/(2r+1)}\ll J_{n}^{\mathrm
{ini}}\ll
n^{2r/(2r+1)}(\log n)^{-1}$ where $r\geq1$. By using the adaptive group
lasso estimator as the initial estimator, Assumption~\ref{Assumption1}
requires that $J_{n}^{\mathrm{ini}}\ll\{n\log(n)\}^{1/2}$. Hence $%
n^{1/(2r+1)}\ll J_{n}^{\mathrm{ini}}\ll\{n\log(n)\}^{1/2}$. We
therefore can
let $J_{n}^{\mathrm{ini}}\asymp n^{(1+\vartheta)/(2r+1)}$, where
$\vartheta$
is any small positive number close to $0$. This increase in the number of
basis functions ensures undersmoothing in the first step in order that the
uniform difference between the two-step and the oracle estimators become
asymptotically negligible. Based on Assumptions \ref{Assumption1} and
\ref%
{Assumption2}, the tuning parameter $\lambda_{n}$ needs to satisfy $%
n^{-1/2}(J_{n}^{\mathrm{ini}})^{1/2}\sqrt{\log(pJ_{n}^{\mathrm{ini}})}%
\{\min_{\ell\in I_{2}}(w_{n\ell})\}^{-1}\ll\lambda_{n}\ll1$.
\end{remark}

\begin{remark}
\label{re:bic2} The number of interior knots has the same order requirement
as the number of basis functions. In the first step, with the undersmoothing
requirement as discussed in Remark~\ref{Remark2}, we let the number of
interior knots $N^{\mathrm{ini}}=\lfloor cn^{( 1+0.01) /( 2q+1) }\rfloor$,
where $c$ is a constant, by assuming that $r=q$. In the simulations, we
let $%
c=2$. In the second-step estimation, we use BIC to select the number of
knots $N^{S}$, so the optimal $N^{S}$ ranges in $[ \lfloor n^{1/( 2q+1)
}\rfloor,\lfloor2n^{1/( 2q+1) }\rfloor] $ by minimizing BIC: BIC$(
N^{S}) =2L_{n}^{\mathcal{S}}( {\widehat{{\bolds{\gamma}}}}%
_{,1}^{\mathcal{S}}) +d( N^{S}+q) $log$( n) $.
\end{remark}

\subsection{Simultaneous confidence bands}\label{sec3.5}

In this section, we propose a SCB for $\alpha_{\ell1}(x_{1})$
by studying the asymptotic behavior of the maximum of the normalized
deviation of the spline functional estimate. To construct asymptotic SCBs
for $\alpha_{\ell1}(x_{1})$ over the interval $x_{1}\in{}[0,1]$ with
confidence level $100(1-\alpha)\%$, $\alpha\in(0,1)$, we need to
find two
functions $l_{\ell n}(x_{1})$ and $u_{\ell n}(x_{1})$ such that
%
\begin{equation}
\lim_{n\rightarrow\infty}P\bigl(l_{\ell n}(x_{1})\leq
\alpha_{\ell
1}(x_{1})\leq u_{\ell n}(x_{1})
\mbox{ for all }x_{1}\in{}[ 0,1]\bigr)=1-\alpha. \label{EQ:ACB}
\end{equation}
In practice, we consider a variant of (\ref{EQ:ACB}) and construct
SCBs over
a subset $S_{n,1}$ of $[0,1]$ with $S_{n,1}$ becoming denser as $%
n\rightarrow\infty$. We, therefore, partition $[0,1]$ according to $L_{n}$
equally spaced intervals based on $0<\xi_{0}<\xi_{1}<\cdots<\xi
_{L_{n}}<\xi_{L_{n}+1}=1$ where $L_{n}\rightarrow\infty$ as
$n\rightarrow
\infty$. Let $S_{n,1}=(\xi_{0},\ldots,\xi_{L_{n}})$. Define
$ 
d_{L_{n}}(\alpha)=1-\{2\log(L_{n}+1)\}^{-1}[\log\{-(1/2)\log
(1-\alpha
)\}+(1/2)\{\log\log(L_{n}+1)+\log(4\pi)\}]$, 
and $Q_{L_{n}}(\alpha)=\{2\log(L_{n}+1)\}^{1/2}d_{L_{n}}(\alpha)$.

\begin{theorem}
\label{THM:bands} Under conditions \textup{(C1)--(C5)} in the \hyperref[app]{Appendix}, and $%
L_{n}\asymp J_{n}^{\mathcal{S}}\asymp n^{1/(2r+1)}$ and
$n^{1/(2r+1)}\ll
J_{n}^{\mathrm{ini}}\ll n^{2r/(2r+1)}\{\log(n)\}^{-1}$, we have
\[
\lim_{n\rightarrow\infty}P \Bigl\{ \sup_{x_{1}\in S_{n,1}}\bigl\llvert
\sigma _{n1}^{-1}(x_{1})\bigl\{\widehat{
\alpha}_{\ell1}^{\mathcal
{S}}(x_{1})-\alpha
_{\ell1}(x_{1})\bigr\}\bigr\rrvert \leq Q_{L_{n}}(
\alpha) \Bigr\} =1-\alpha,
\]
and thus an asymptotic $100(1-\alpha)\%$ confidence band for $\alpha
_{\ell
1}(x_{1})$ over $x_{1}\in S_{n,1}$ is
%
\begin{equation}
\widehat{\alpha}_{\ell1}^{\mathcal{S}}(x_{1})\pm\sigma
_{n1}(x_{1})Q_{L_{n}}(\alpha). \label{EQ:CB}
\end{equation}
\end{theorem}

\begin{remark}
Compared to the pointwise confidence intervals with width
$2Z_{1-\alpha
/2}\sigma_{n}(x_{1})$, the width of the confidence bands (\ref
{EQ:CB}) is
inflated by a rate $\{2\log(L_{n}+1)\}^{1/2}d_{L_{n}}(\alpha
)/Z_{1-\alpha
/2}$, where $Z_{1-\alpha/2}$ is the cut-off point of the $100(1-\alpha
)$th percentile of the standard normal.
\end{remark}

\subsection{Bootstrap smoothing for calculating the standard
error}\label{sec3.6}

Theorem~\ref{THM:bands} establishes a thresholding value
$Q_{L_{n}}(\alpha)$ for the SCB. One critical question is how to
estimate the standard
deviation $\sigma_{n1}(x_{1})$ in order to construct the SCB. We can
use a sample estimate of $\sigma_{n1}(x_{1})$ according to the asymptotic
formula given in~(\ref{DEF:sig2nx1}), which may have approximation
error and thus lead to inaccurate results for inference. The bootstrap
estimate of
the standard deviation provides an alternative way. We here propose a
bootstrap smoothed confidence band by adopting the nonparametric bootstrap
smoothing idea from \citet{efron:2013}, which can eliminates
discontinuities in jumpy estimates. The procedure is described as follows.

Let $\mathbf{D}=\{\mathbf{D}_{1},\ldots,\mathbf{D}_{n}\}$
be the
data we have, where $\mathbf{D}_{i}=\{Y_{i},\mathbf{X}_{i},(T_{i\ell
},\ell
\in\widehat{I}_{1})\}$. Denote $\mathbf{D}^{\ast}=\{
\mathbf
{D}%
_{1}^{\ast},\ldots,\mathbf{D}_{n}^{\ast}\}$ as a nonparametric bootstrap
sample from $\{\mathbf{D}_{1},\ldots,\mathbf{D}_{n}\}$, and $\mathbf{D}
_{(j)}^{\ast}=\{\mathbf{D}_{(j)1}^{\ast},
\ldots,\mathbf{D}_{(j)n}^{\ast}\}$ as the $j$th bootstrap sample
in $B$ draws. Let $\widehat{\alpha}_{\ell1,(j)}^{\mathcal{\ast S}}(x_{1})$
be the two-step estimator of $\alpha_{\ell1}(x_{1})$ by using the
data $%
\mathbf{D}_{(j)}^{\ast}$. We first present an empirical standard deviation
by the traditional resampling method which is given as
%
\begin{equation}
\widehat{\sigma}_{\ell1,B}(x_{1})=\Biggl[\sum
_{j=1}^{B}\bigl\{\widehat{%
\alpha}_{\ell1,(j)}^{\mathcal{\ast S}}(x_{1})-\widehat{\alpha
}_{\ell
1,\cdot}^{\mathcal{\ast S}}(x_{1})\bigr\}^{2}\Big/(B-1)
\Biggr]^{1/2}, \label{EQ:SE2}
\end{equation}
where $\widehat{\alpha}_{\ell1,\cdot}^{\mathcal{\ast S}%
}(x_{1})=\sum_{j=1}^{B}\widehat{\alpha}_{\ell1,(j)}^{\mathcal
{%
\ast S}}(x_{1})/B$. Then a $100(1-\alpha)\%$ unsmoothed bootstrap SCB
for $%
\alpha_{\ell1}(x_{1})$ over $x_{1}\in S_{n,1}$ is given as
%
\begin{equation}
\widehat{\alpha}_{\ell1}^{\mathcal{S}}(x_{1})\pm\widehat{
\sigma }_{\ell
1,B}(x_{1})Q_{L_{n}}(\alpha).
\label{EQ:CB2}
\end{equation}
Another choice is the smoothed bootstrap SCB which eliminates
discontinuities in the estimates [\citet{efron:2013}]. Let
\[
\widetilde{\alpha}_{\ell1}^{\mathcal{S}}(x_{1})=\sum
_{j=1}^{B}%
\widehat{
\alpha}_{\ell1,(j)}^{\mathcal{\ast S}}(x_{1})/B
\]
be the smoothed estimate of $\alpha_{\ell1}(x_{1})$ obtained by averaging
over the bootstrap replications. Let $C_{(j)i}^{\ast}=\#\{\mathbf{D}%
_{(j)i^{\prime}}^{\ast}=\mathbf{D}_{i}\}$ be the number of elements
in $%
\mathbf{D}_{(j)i^{\prime}}^{\ast}$ equaling $\mathbf{D}_{i}$.

\begin{proposition}
\label{PROP:bootstrap} At each point $x_{1}\in S_{n,1}$, the nonparametric
delta-\break method estimate of the standard deviation for the smoothed bootstrap
statistic $\widetilde{\alpha}_{\ell1}^{\mathcal{S}}(x_{1})$ is
$\widetilde{%
\sigma}_{\ell1}(x_{1})=\{\sum_{i=1}^{n}\operatorname{cov}_{i}^{2}(x_{1})\}
^{1/2}$%
, where $\operatorname{cov}_{i}(x_{1})=\operatorname{cov}_{\ast}\{C_{(j)i}^{\ast},
\break \widehat{\alpha
}_{\ell
1,(j)}^{\mathcal{\ast S}}(x_{1})\}$ which is the bootstrap covariance
between $%
C_{(j)i}^{\ast}$ and $\widehat{\alpha}_{\ell1,(j)}^{\mathcal{\ast
S}%
}(x_{1})$.
\end{proposition}

The proof of Proposition~\ref{PROP:bootstrap} essentially follows the same
arguments as the proof for Theorem~1 in \citet{efron:2013}. Based on
Proposition~\ref{PROP:bootstrap}, to construct the smoothed
bootstrap SCB, we use the nonparametric estimate of the standard deviation
given as
%
\begin{equation}
\widetilde{\sigma}_{\ell1,B}(x_{1})=\Biggl\{\sum
_{i=1}^{n}\widehat {\operatorname{cov}}%
_{i,B}^{2}(x_{1})
\Biggr\}^{1/2}, \label{EQ:SE3}
\end{equation}
where
\[
\widehat{\operatorname{cov}}_{\ell i,B}(x_{1})=\sum
_{j=1}^{B}\bigl(C_{(j)i}^{\ast
}-C_{\cdot i}^{\ast}
\bigr) \bigl(\widehat{\alpha}_{\ell1,(j)}^{\mathcal{\ast
S}%
}(x_{1})-
\widehat{\alpha}_{\ell1,\cdot}^{\mathcal{\ast S}}(x_{1})\bigr)/B
\]
with $C_{\cdot i}^{\ast}=\sum_{j=1}^{B}C_{(j)i}^{\ast}/B$.
The $%
100(1-\alpha)\%$ smoothed bootstrap SCB for $\alpha_{\ell1}(x_{1})$
over $%
x_{1}\in S_{n,1}$ is given as
%
\begin{equation}
\widetilde{\alpha}_{\ell1}^{\mathcal{S}}(x_{1})\pm
\widetilde {\sigma} 
_{\ell1,B}(x_{1})Q_{L_{n}}(
\alpha). \label{EQ:CB3}
\end{equation}

\section{A simulation study}\label{sec7}

In this section, we present a simulation study to evaluate the finite sample
performance of our proposed penalized estimation procedure and the
simultaneous confidence bands. More numerical studies are located in the
supplementary materials [\citet{ma.carroll.liang.xu2015}].

\begin{example}\label{ex1}
In this example, we use $1286$ SNPs located on the sixth
chromosome from the Framingham Heart Study to simulate the binary response
from the logistic model
%
\begin{equation}\quad
\operatorname{logit}\bigl\{P(Y_{i}=1|\mathbf{X}_{i},
\mathbf{T}_{i})\bigr\}=\sum_{\ell
=1}^{p}
\alpha_{\ell}(\mathbf{X}_{i})T_{i\ell}=\sum
_{\ell
=1}^{p}\Biggl\{\alpha_{\ell0}+\sum
_{k=1}^{2}\alpha_{\ell
k}(X_{ik})
\Biggr\}T_{i\ell}, \label{MOD:LOGIT}
\end{equation}
with the four SNPs ss66063578, ss66236230, ss66194604 and ss66533844
selected from the real data analysis in Section~\ref{sec8} as important
covariates and the other SNPs as unimportant covariates, so that $s=4$ (the
number of important covariates), $p=1286$ and the sample size $n=300$. The
three possible allele combinations are coded as 1, 0 and $-1$ for each SNP.
The covariates $X_{ik}$, $k=1,2$, are simulated environmental effects, which
are generated from independent uniform distributions on $[0,1]$. We generate
the coefficient functions as $\alpha_{10}=0.5$, $\alpha
_{11}(x_{1})=4\cos
(2\pi x_{1})$, $\alpha_{12}(x_{2})=5\{(2x_{2}-1)^{2}-1/3\}$, $\alpha
_{20}=0.5$, $\alpha_{21}(x_{1})=6x_{1}-3$, $\alpha_{22}(x_{2})=4\{
\sin
(2\pi x_{2})+\cos(2\pi x_{2})\}$, $\alpha_{30}=0.5$, $\alpha
_{31}(x_{1})=4\sin(2\pi x_{1})$, $\alpha_{32}(x_{2})=6x_{2}-3$,
$\alpha
_{40}=0.5$, $\alpha_{41}(x_{1})=4\cos(2\pi x_{1})$, $\alpha
_{42}(x_{2})=5\{(2x_{2}-1)^{2}-1/3\}$ and $\alpha_{\ell}(\mathbf{X}_{i})=0
$ for $l=5,\ldots,1286$. We conducted $500$ replications for each
simulation. We fit the data with the GACM (\ref{MOD:LOGIT}) by using the
adaptive group lasso (AGL) and group lasso (GL). In the literature, the
generalized varying coefficient model [GVCM; \citet{lian:2012}], which
considers one index variable in the coefficient function for each
predictor $%
T_{i\ell}$, has been widely used to study nonlinear interactions. To
apply the GVCM method [\citet{lian:2012}] in this setting, we first perform
principal component analysis (PCA) on $\mathbf{X}_{i}$ and then use the
first principal component as the index variable in the GVCM. Then we apply
the AGL and GL methods to the GVCM: $\operatorname{logit}\{P(Y_{i}=1|\mathbf
{X}_{i},\mathbf{%
T}_{i})\}=\sum_{\ell=1}^{p}\alpha_{\ell}(U_{i})T_{i\ell}$,
where $U_{i}$ is the first principal component obtained by PCA on
$\mathbf{X}%
_{i}$. Moreover, we also fit the data with the parametric logistic
regression by assuming linear coefficient functions (\ref{MOD:Linear}) with
the AGL method. We also compare our proposed method with the conventional
screening method by parametric logistic regression for Genome-Wide
Association Studies [GWAS; \citet{Murcray.Lewinger.James.Gauderman:2008}].
In the screening method, we fit a logistic model for each SNP: $\operatorname{logit}%
\{P(Y_{i}=1|\mathbf{X}_{i},T_{i\ell})\}=\alpha_{0}+\bolds{\alpha}^{\mathrm{T}}\mathbf{X}_{i}+\beta_{\ell}T_{i\ell
}+\sum_{k=1}^{2}\beta_{\ell k}X_{ik}T_{i\ell}$, for $\ell
=1,\ldots,1286$. Then we conduct a likelihood ratio test for the genetic
and interaction effects of $H_{0}:\beta_{\ell}=\beta_{\ell1}=\beta
_{\ell2}=\beta_{\ell3}=0$. Let $\alpha_{0}=0.05$ be the overall
type I
error for the study and $M=1286$ be the number of SNPs in this study. We
apply the multiple testing correction procedure for GWAS with $H_{0}$
rejected when the $p$-value${}<\alpha_{0}/M_{\mathrm{eff}}$, where
$M_{\mathrm{eff}}
$ is the Cheverud--Nyholt estimate of the effective number of tests %
[\citet{Cheverud:2001,Nyholt:2004}] calculated by $M_{\mathrm{eff}%
}=1+M^{-1}\sum_{j=1}^{M}\sum_{k=1}^{M}(
1-r_{jk}^{2})$ and $r_{jk}$ are the correlation coefficients of the
SNPs, and we obtain $M_{\mathrm{eff}}=1275.65$.

\begin{table}
\caption{Variable selection and estimation results by
the adaptive group lasso and the group lasso with the GACM and GVCM,
respectively, and parametric logistic regression with adaptive group lasso
and screening methods based on $500$ replications. The columns of C, O
and I
show the percentage of correct-fitting, over-fitting and incorrect-fitting.
The columns TP, FP and MR show true positives, false positives and model
errors, respectively}\label{TAB:vc2}
\begin{tabular*}{\textwidth}{@{\extracolsep{\fill}}lccccccc@{}}
\hline
& & \textbf{C} & \textbf{O} & \textbf{I} & \textbf{TP} & \textbf{FP} & \textbf{MR} \\
\hline
GACM & AGL & $0.410$ & $0.460$ & $0.130$ & $3.860$ & \phantom{0}$0.870$ & $0.059$
\\
& GL & $0.140$ & $0.764$ & $0.096$ & $3.904$ & $\phantom{0}2.540$ & $0.083$ \\[3pt]
GVCM & AGL & $0.030$ & $0.000$ & $0.970$ & $1.636$ & \phantom{0}$5.685$ & $0.142$
\\
& GL & $0.060$ & $0.000$ & $0.940$ & $2.076$ & $20.670$ & $0.120$ \\[3pt]
Logistic regression & AGL & $0.000$ & $0.000$ & $1.000$ & $1.872$ & \phantom{0}$1.174$
& $0.159$ \\
& Screening & $0.000$ & $0.000$ & $1.000$ & $1.056$ & \phantom{0}$0.786$ & $0.141$
\\
\hline
\end{tabular*}
\end{table}

Table~\ref{TAB:vc2} presents the percentages of correct-fitting (C) (exactly
the important covariates are selected), over-fitting (O) (both the important
covariates and some unimportant covariates are selected) and
incorrect-fitting (I) (some of the important covariates are not selected),
the average true positives (TP), that is, the average number of selected
covariates among the important covariates, the average false positives (FP),
that is, the average number of selected covariates among the unimportant
covariates, and the average model errors (MR), the latter defined as $%
\sum_{i=1}^{n}\{\widehat{\mu}_{i}(\mathbf{X}_{i},\mathbf{T}%
_{i})-\mu_{i}(\mathbf{X}_{i},\mathbf{T}_{i})\}^{2}/n$, where
$\widehat
{\mu}%
_{i}(\mathbf{X}_{i},\mathbf{T}_{i})$ and $\mu_{i}(\mathbf
{X}_{i},\mathbf
{T}%
_{i})$ are the estimated and true conditional means for $Y_{i}$,
respectively. We see that by fitting the proposed GACM, the GL method has
larger percentage of over-fitting as well as larger average false positives
than the AGL methods. The AGL improves the correct-fitting percentage
by $%
26\%$. As a result, the AGL reduces the model fitting error by $%
(0.083-0.059)/0.059=40.7\%$ compared to the GL method. Moreover, both the
logistic model and the GVCM fail to identify those important covariates with
incorrect-fitting percentage close to or being $1$. Furthermore, by using
the screening method with logistic regression, the average true
positive is $%
1.056$, which is much less than $4$ (the number of those important SNPs).
This further illustrates that the traditional screening method is not an
effective tool to identify important genetic factors in this context. In
addition, we observe that the results for the AGL method in Table~\ref%
{TAB:vc2} are comparable to the results in Table~S.1 of
Example~2 (in the supplementary materials)
at $p=1000$ with the simulated SNPs in terms of having similar
correct-fitting percentages and MR values.

Next, we investigate the empirical coverage rates of the unsmoothed and
smoothed SCBs given in (\ref{EQ:CB2}) and (\ref{EQ:CB3}). To
calculate the
unsmoothed and smoothed bootstrap standard deviations (\ref{EQ:SE2})
and (%
\ref{EQ:SE3}), we use $B=500$ bootstrap replications. The confidence bands
are constructed at $L_{n}=20$ equally spaced points. At $95\%$ confidence
level, Table~\ref{TAB:boot2} reports the empirical coverage rates (cov)
and the sample averages of median and mean standard deviations
(sd.median and sd.mean), respectively, for the unsmoothed SCB (\ref%
{EQ:CB2}) and smoothed SCB (\ref{EQ:CB3}) for coefficient functions
$\alpha
_{\ell1}(x_{1})$, $\ell=1,2,3,4$. We see that the smoothed bootstrap
method leads to better performance, having empirical coverage rates closer
to the nominal confidence level $0.95$.
\end{example}

\begin{table}
\caption{
The empirical coverage rates (cov) and
the sample average of median and mean of the standard deviations
(sd.median and sd.mean) for the unsmoothed SCB (\protect\ref{EQ:CB2}%
) and smoothed SCB (\protect\ref{EQ:CB3}) for the coefficient functions
$%
\protect\alpha_{\ell1}(x_{1})$ for $\ell=1,2,3,4$}\label{TAB:boot2}
\begin{tabular*}{\textwidth}{@{\extracolsep{\fill}}lcccccc@{}}
\hline
& \multicolumn{3}{c}{\textbf{Unsmoothed bootstrap}} & \multicolumn{3}{c}{\textbf{Smoothed
bootstrap}} \\ [-6pt]
& \multicolumn{3}{c}{\hrulefill} & \multicolumn{3}{c}{\hrulefill} \\
& \textbf{cov} & \textbf{sd.median} & \textbf{sd.mean} & \textbf{cov} &
\textbf{sd.median} &
\multicolumn{1}{c@{}}{\textbf{sd.mean}} \\
\hline
$\alpha_{11}$ & 0.610 & 0.689 & 0.809 & 0.818 & 0.735 & 0.982 \\
$\alpha_{21}$ & 0.628 & 0.563 & 0.725 & 0.846 & 0.666
& 0.932 \\
$\alpha_{31}$ & 0.636 & 0.736 & 0.832 & 0.869 & 0.837 & 1.053 \\
$\alpha_{41}$ & 0.646 & 0.768 & 0.843 & 0.882 & 0.891
& 1.064 \\
\hline
\end{tabular*}
\end{table}

\section{Data application}\label{sec8}

We illustrate our method via analysis of the Framingham Heart Study %
[\citet{Dawber.Meadors.Moore:1951}] to investigate the effects of
G${}\times{}$E
interactions on obesity. People are defined as obese when their body mass
index (BMI) is 30 or greater: this is the definition of being obese
made by
the U.S. Centers for Disease Control and Prevention; see
\url{http://www.cdc.gov/obesity/adult/defining.html}. We defined the response
variable to be $Y=1$ for BMI${}\geq30$; and $Y=0$ for BMI${}<30$. We use $%
X_{1}=$ sleeping hours per day; $X_{2}=$ activity hours per day; and $X_{3}=$
diastolic blood pressure as the environmental factors, and use single
nucleotide polymorphisms (SNPs) located in the sixth chromosome as the
genetic factors. The three possible allele combinations are coded as 1, 0
and $-1$. As in the simulation, we thus are fitting a multiplicative risk
model in the SNPs. For details on genotyping, see
\url{http://www.ncbi.nlm.nih.gov/projects/gap/cgi-bin/study.cgi?studyid=phs000007.v3.p2}.
A total of
$1286$ SNPs remain in our analysis after eliminating SNPs with minor allele
frequency $<$0.05, those with departure from Hardy--Weinberg
equilibrium and
those having correlation coefficient with the response between $-0.1$
and $%
0.1$. We have $n=300$ individuals left in our study after deleting
observations with missing values.

To see possible nonlinear main effects of the environmental factors, we
first fit a generalized additive model by using $X_{1}$, $X_{2}$ and $X_{3}$
as predictors such that
%
\begin{equation}
E(Y_{i}|\mathbf{X}_{i}\mathbf{,T}_{i})=g^{-1}
\bigl\{\eta(\mathbf {X}_{i})%
\bigr\}\qquad\mbox{with }\eta(
\mathbf{X}_{i})=m_{0}+\sum%
_{k=1}^{3}m_{k}(X_{ik}).
\label{MOD:GAM}
\end{equation}
Figure~S.1 given in the supplementary material [\citet
{ma.carroll.liang.xu2015}] depicts the plots of $\widehat{m}_{k}(
\cdot
) $ for $k=1,2,3$ by one-step cubic spline estimation. Clearly the
estimate of each nonparametric function has a nonlinear pattern. We
refer to Section~S.2 for the detailed description of this figure.
Based on the plots shown in Figure~S.1, we fit the GACM model
%
\begin{equation}
\eta(\mathbf{X}_{i},\mathbf{T}_{i})=\sum
_{\ell=1}^{1287}\Biggl\{ \alpha _{\ell0}+\sum
_{k=1}^{3}\alpha_{\ell k}(X_{ik})
\Biggr\}T_{i\ell}, \label{MOD:real}
\end{equation}
where $\mathbf{T}_{i}=(T_{i1},T_{i2},\ldots,T_{i1287})^{\mathrm{T}}$
with $%
T_{i1}=1$, and $T_{i\ell}$ are the SNP covariates for $\ell=2,\ldots
,1287$. The nonparametric function $\alpha_{\ell k}(\cdot)$ is
estimated by
cubic splines, and the number of interior knots for each step is selected
based on the criterion described in Section~\ref{sec:computation}. We
select variables in model (\ref{MOD:real}) by the proposed adaptive group
lasso (AGL) and the group lasso (GL). To compare the proposed model with
linear models, we perform the group lasso by assuming linear interaction
effects (Linear) such that $\alpha_{\ell}(\mathbf{X}_{i})=\alpha
_{\ell
0}+\sum_{k=1}^{3}\beta_{\ell k}X_{ik}$, and we also perform the
lasso by assuming no interaction effects (No interaction) such that
$\alpha
_{\ell}(\mathbf{X}_{i})=\alpha_{\ell0}$. We also apply the screening
method with parametric logistic regression (Screening) as described in
Example~2. Table~\ref{TAB:NAT2} reports the variable
selection results in these five scenarios. After model selection, we
calculate the estimated leave-one-out cross-validation prediction error
(CVPE) for the model with the selected variables as shown in the last
row of
Table~\ref{TAB:NAT2}. Among the selected SNPs by the AGL method, two SNPs,
rs4714924 and rs6543930, have been scientifically confirmed by \citet%
{Randall.Winkler:2013} to have strong associations with obesity. Moreover,
compared to the linear, no interaction and screening methods, our proposed
AGL with GACM method enables us to identify more genetic factors, which may
be important to the response but missed out by other methods. As a result,
it has the smallest CVPE ($0.078$), so that it significantly improves model
prediction compared to other methods. We also see that the logistic model
that completely ignores interactions has the largest CVPE $(0.152)$. The
screening method has the second largest CVPE $(0.149)$, which is larger than
that of the penalization method ($0.124$) obtained by fitting the same
logistic regression model but including interaction considered. This result
demonstrates that the screening method is not as effective as the
penalization method for analysis of this data set, a result which also
agrees with our simulations.

\begin{table}
\caption{Variable selection results for the group lasso (GL) and the
adaptive group lasso (AGL) in model~(\protect\ref{MOD:real}), the group
lasso by assuming linear interaction effects (linear), the lasso by assuming
no interaction effects (no interaction) and the screening method
(screening). The symbol $\surd$ indicates that the SNP was selected into
the model. The last row shows the cross validation prediction errors (CVPE)}
\label{TAB:NAT2}
\begin{tabular*}{\textwidth}{@{\extracolsep{\fill}}lccccc@{}}
\hline
\textbf{SNPs} & \textbf{GL} & \textbf{AGL} & \textbf{Linear} &
\textbf{No interaction} & \textbf{Screening} \\
\hline
rs9296244 & $\surd$ & $\surd$ & & & \\
rs6910353 & $\surd$ & $\surd$ & & & \\
rs3130813 & $\surd$ & $\surd$ & & & \\
rs9353447 & & & & $\surd$ & $\surd$ \\
rs4714924 & $\surd$ & $\surd$ & $\surd$ & & \\
rs242263 & $\surd$ & $\surd$ & $\surd$ & & $\surd$ \\
rs282123 & $\surd$ & & & & \\
rs282128 & $\surd$ & $\surd$ & & & \\
rs6929006 & & & $\surd$ & & \\
rs9353711 & $\surd$ & & & & \\
rs12199154 & $\surd$ & $\surd$ & & & \\
rs2277114 & $\surd$ & & & & \\
rs749517 & $\surd$ & & & & \\
rs729888 & $\surd$ & & & & \\
rs203139 & & & $\surd$ & & \\
rs6914589 & $\surd$ & $\surd$ & & & \\
rs6543930 & $\surd$ & $\surd$ & & & \\[6pt]
CVPE & $0.099$ & $0.078$ & $0.124$ & $0.152$ & $0.149$ \\
\hline
\end{tabular*}
\end{table}

\begin{figure}

\includegraphics{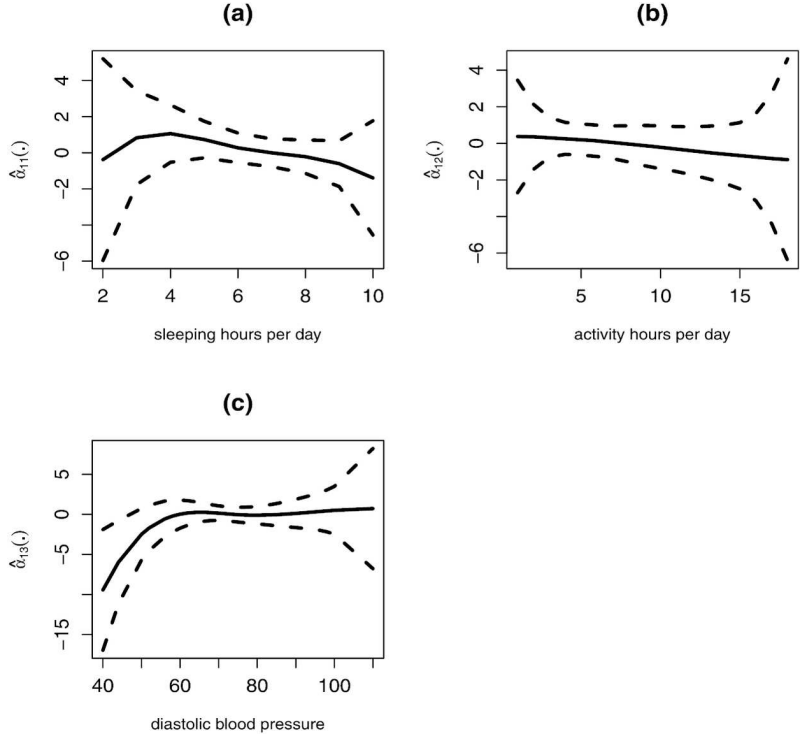}

\caption{Plots of the smoothed two-step estimated functions $\protect%
\widetilde{\protect\alpha}_{1k}^{\mathcal{S}}(\cdot)$ for
$k=1,2,3$ and
the associated $95\%$ SCBs based on model \protect\eqref{eq:finalMod}.}
\label{FIG:real}
\end{figure}

Next we fit the final GACM selected variables from the AGL procedure as
%
\begin{equation}
\eta(\mathbf{X}_{i},\mathbf{T}_{i})=\sum
_{\ell=1}^{10}\Biggl\{ \alpha _{\ell0}+\sum
_{k=1}^{3}\alpha_{\ell k}(X_{ik})
\Biggr\}T_{i\ell}. \label{eq:finalMod}
\end{equation}
To illustrate the main effects of the environmental factors, Figure~\ref
{FIG:real} plots the smoothed two-step estimated functions $\widetilde{
\alpha}_{1k}^{\mathcal{S}}(\cdot)$ of the functions $\alpha_{1k}^{%
\mathcal{S}}(\cdot)$, for $k=1,2,3$, and the associated $95\%$ smoothed
SCBs (upper and lower solid lines). The plots of the functional estimates
have the same nonlinear change patterns as the corresponding plots in
Figure~S.1, although because of the addition of the SCBs, the scale
of the plot has changed.

\begin{figure}

\includegraphics{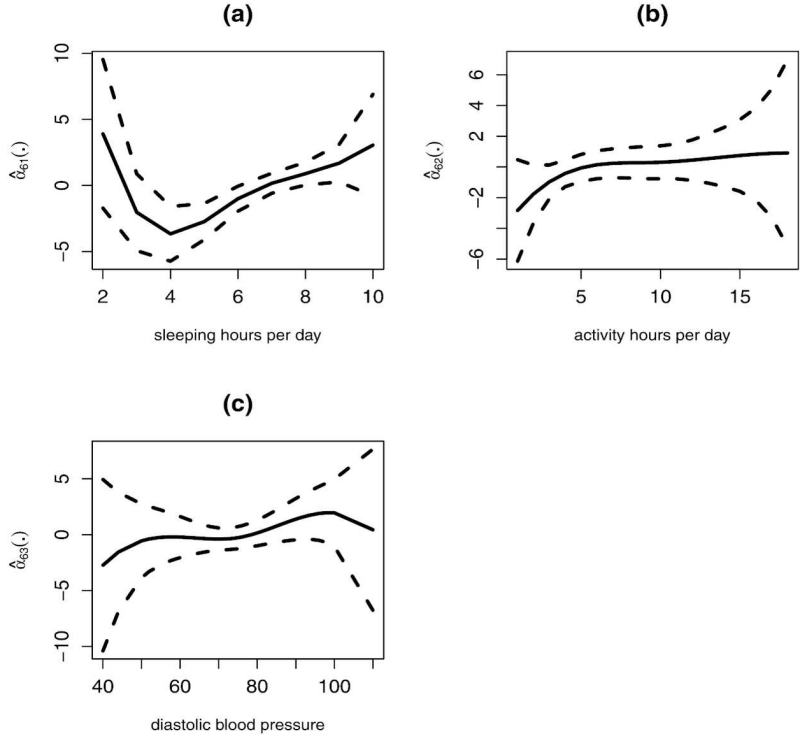}

\caption{Plots of the smoothed two-step estimated functions $\protect%
\widetilde{\protect\alpha}_{5k}^{\mathcal{S}}(\cdot)$ for
$k=1,2,3$ and
the associated $95\%$ SCBs based on model \protect\eqref{eq:finalMod}.}
\label{FIG:real2}
\end{figure}

\begin{figure}

\includegraphics{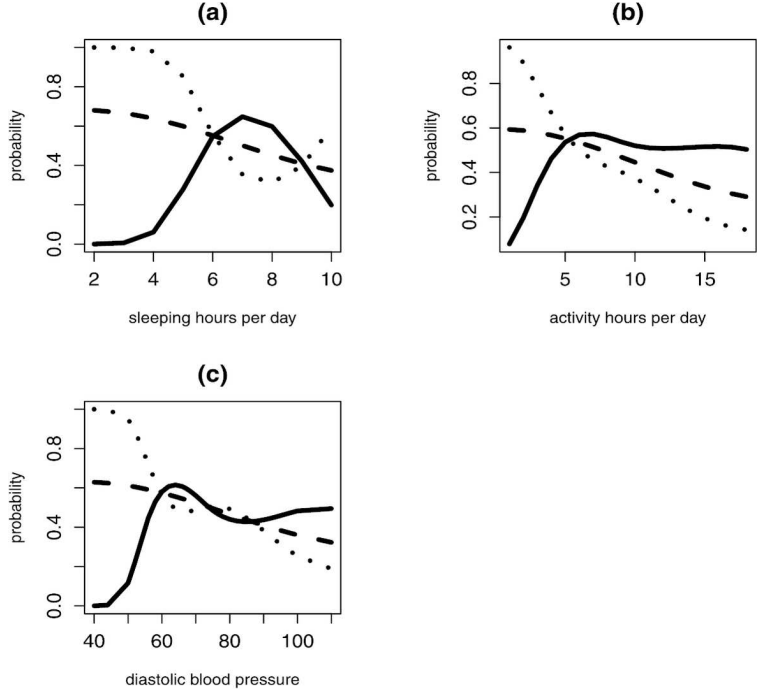}

\caption{Plots of the estimated conditional probability of obesity against
each environmental factor by letting $T_{i\ell}=0$ for $\ell\neq5$. With
A being the minor allele, the curves are aa (solid line), Aa (dashed line)
and AA (dotted line), based on model \protect\eqref{eq:finalMod}.}
\label{FIG:real3}
\end{figure}

To illustrate the effects of the genetic factors changing with the
environmental factors, in Figure~\ref{FIG:real2} we plot the smoothed
two-step estimated functions $\widetilde{\alpha}_{6k}^{\mathcal
{S}}(\cdot)$
and the associated $95\%$ smoothed SCBs of the coefficient functions
$\alpha
_{6k}^{\mathcal{S}}(\cdot)$ for the SNP rs242263. To further demonstrate
how the probability of developing obesity changes with the environmental
factors for each category of SNP rs242263, Figure~\ref{FIG:real3}
plots the
estimated conditional probability of obesity against each environmental
factor by letting $T_{i\ell}=0$ for $\ell\neq6$. Letting A be the minor
allele, the curves are for aa (solid line), Aa (dashed line) and AA (dotted
line). Figure~\ref{FIG:real2} indicates different changing patterns of the
interaction effects under different environments. For example, sleeping
hours seem to have an overall more significant interaction effect with this
particular SNP than the other two variables. The effect of this SNP changes
from positive to negative and then to positive again as the sleeping hours
increase. The coefficient functions of the SNP have an increasing pattern
along with the activity hours and diastolic blood pressure, respectively.
From Figure~\ref{FIG:real3}, we observe that there are stronger differences
among the levels AA, Aa, and aa of SNP rs242263 for both large and small
values of the environmental factors.There are other interesting results
worth further study. For example, in the 2--6 hours per day sleeping range,
the AA group (dotted lines) have much higher rates of obesity than the aa
group (solid line), but the opposite occurs in the 6--9 hour range. For those
with low amounts of activity per day, again the AA group is more obese than
the aa group, while when activity increases, the AA group is less obese than
the aa group. A similar noticeable difference occurs between the $<$60
diastolic blood pressure group, those who are hypotensive, and the $>$90
group, those who are hypertensive, although there are few subjects in the
former group.

\section{Discussions}\label{sec:disc}

The generalized additive coefficient model (GACM) proposed by \citet%
{xue.yang:2006} and \citet{xue.liang:2010} has been demonstrated to be a
powerful tool for studying nonlinear interaction effects of variables. To
promote the use of the GACM in modern data applications such as
gene-environment (G${}\times{}$E) interaction effects in GWAS, we have proposed
estimation and inference procedures for the GACM when the dimension of the
variables is high. Specifically, we have devised a groupwise penalization
method in the GACM for simultaneous model selection and estimation. We
showed by numerical studies that we can effectively identify important
genetic factors by using the proposed nonparametric model while traditional
generalized parametric models such as logistic regression model fails
to do
so when nonlinear interactions exist. Moreover, by comparing with the
conventional screening method with logistic regression as commonly used in
the GWAS community, our proposed groupwise penalization method with the GACM
has been demonstrated to be more effective for variable selection and model
estimation. After identifying those important covariates, we have further
constructed simultaneous confidence bands for the nonzero coefficient
functions based on a refined two-step estimator. We estimate the standard
deviation of the functional estimator by a smoothed bootstrap method as
proposed in \citet{efron:2013}. The method was shown to have good numerical
performance by reducing variability as well as improving the empirical
coverage rate of the proposed simultaneous confidence bands. Our
methods can
be extended to longitudinal data settings through marginal models or
mixed-effects models. More work, however, is needed to understand the
properties of the estimators in such new settings. Moreover, extending this
work to the setting with the dimensions for both genetic and environmental
factors growing with the sample size can be a future project to be
considered. Some associated theoretical properties with respect to model
selection and estimation as well as inference need to be carefully
investigated.

\begin{appendix}
\section*{Appendix}\label{app}

Denote the space of the $q$th order smooth functions as $%
C^{(q)}([0,1])=\{\phi|\phi^{(q)}\in C[0,1]\}$. For any $s\times s$
symmetric matrix $\mathbf{A}$, denote its $L_{q}$ norm as $\Vert
\mathbf
{A}%
\Vert_{q}=\max_{\varsigma\in R^{s},\llVert \varsigma\rrVert
_{2}=1}\Vert\mathbf{A}\varsigma\Vert_{q}$. Let $\Vert\mathbf
{A}\Vert
_{\infty}=\max_{1\leq i\leq s}\sum_{j=1}^{s}|a_{ij}|$. For a
vector $\mathbf{a}$, let $\Vert\mathbf{a}\Vert_{\infty}=\max_{1\leq
i\leq
s}|a_{i}|$.

Let $C^{0,1}(\mathcal{X}_{w})$ be the space of Lipschitz continuous
functions on $\mathcal{X}_{w}$, that is,
\[
C^{0,1}(\mathcal{X}_{w})=\biggl\{\varphi:\Vert\varphi
\Vert_{0,1}=\sup_{w\neq
w^{\prime},w,w^{\prime}\in\mathcal{X}_{w}}\frac{|\varphi
(w)-\varphi
(w^{\prime})|}{|w-w^{\prime}|}<+\infty\biggr
\},
\]
in which $\Vert\varphi\Vert_{0,1}$ is the $C^{0,1}$-norm of $\varphi$.
Denote $q_{j}(\eta,y)=\partial^{j}Q\{g^{-1}(\eta),y\}/\partial\eta
^{j}$, so that
\begin{eqnarray*}
q_{1}(\eta,y) &=&\frac{\partial}{\partial\eta}Q\bigl\{g^{-1}(\eta ),y
\bigr\}=-\bigl\{y-g^{-1}(\eta)\bigr\}\rho_{1}(\eta),
\\
q_{2}(\eta,y) &=&\frac{\partial^{2}}{\partial\eta^{2}}Q\bigl\{ g^{-1}(\eta ),y
\bigr\}=\rho_{2}(\eta)-\bigl\{y-g^{-1}(\eta)\bigr\}
\rho_{1}^{\prime}(\eta),
\end{eqnarray*}
where $\rho_{j}(\eta)=\{\dot{g}^{-1}(\eta)\}^{j}/V\{g^{-1}(\eta)\}$.

\subsection{Assumptions}\label{appa1}

Throughout the paper, we assume the following regularity
conditions:
\begin{longlist}[(C5)]
\item[(C1)] The joint density of $\mathbf{X}$, denoted by $f(\mathbf{x)}$,
is absolutely continuous, and there exist constants $0<c_{f}\leq
C_{f}<\infty
$, such that $c_{f}\leq$min$_{\mathbf{x}\in[0,1]^{d}}f(\mathbf
{x)}\leq \max_{\mathbf{x}\in[0,1]^{d}}f(\mathbf{x)}\leq C_{f}$.\vspace*{1pt}

\item[(C2)] The function $V$ is twice continuously differentiable, and the
link function $g$ is three times continuously differentiable. The
function $%
q_{2}(\eta,y)<0$ for $\eta\in R$ and $y$ in the range of the response
variable.

\item[(C3)] For $1\leq\ell\leq p$, $1\leq k\leq d$, $\alpha_{\ell
k}^{(r-1)}(x_{k})\in C^{0,1}[0,1]$, for given integer $r\geq1$. The
spline order satisfies $q\geq r$.

\item[(C4)] Let $\varepsilon_{i}=Y_{i}-\mu( \mathbf{X}_{i},\mathbf{T}_{i})
,1\leq i\leq n$. The random variables $\varepsilon_{1},\ldots
,\varepsilon
_{n}$ are i.i.d.
with $E( \varepsilon_{i})
=0$ and $\operatorname{var}( \varepsilon_{i}\vert\mathbf{X}_{i},\mathbf{T}_{i})
=\sigma
^{2}( \mathbf{X}_{i},\mathbf{T}_{i}) $. Furthermore, their tail
probabilities satisfy $P( \vert\varepsilon_{i}\vert>x) <K\exp(
-Cx^{2}) $%
, $i=1,\ldots,n$, for all $x\geq0$ and for some positive constants
$C$ and
$K$.

\item[(C5)] The eigenvalues of $E ( \mathbf{T}_{I_{1}}\mathbf{T}
_{I_{1}}^{\mathrm{T}}\vert\mathbf{X}=\mathbf{x}  ) $, where
$\mathbf{%
T}_{I_{1}}= ( T_{\ell},\ell\in I_{1} ) ^{\mathrm{T}}$, are
uniformly bounded away from $0$ and $\infty$ for all $\mathbf{x}\in%
[0,1]^{d}$. There exist constants $0<c_{1}<C_{1}<\infty$, such that $%
c_{1}\leq E ( T_{\ell}^{2}\vert\mathbf{X}=\mathbf{x}
 )
\leq
C_{1}$, for all $\mathbf{x}\in[0,1]^{d}$, $\ell\in I_{2}$.
\end{longlist}

Conditions (C1)--(C5) are standard conditions for nonparametric estimation.
Condition (C1) is the same as condition (C1) in \citet{xue.yang:2006} and
condition (C5) in \citet{xue.liang:2010}. The first condition in (C2) gives
the assumptions on $V$ and the link function $g$, which can be found in
condition (E) of \citet{lam.fan:2008}. The second condition in (C2)
guarantees that the negative quasi-likelihood function $Q\{g^{-1}(\eta
),y\}$
is convex in $\eta\in R$, which is also given in condition (D) of
\citet
{lam.fan:2008} and (a) of condition 1 in \citet
{carroll.fan.gijbels.wand:1997}%
. Condition (C3) is typical for polynomial spline smoothing; see the same
condition given in Section~5.2 of \citet{huang:2003}. Condition (C4) is the
same as assumption (A2) given in \citet{Huang.Horowitz.Wei:2010}. Condition
(C5) is given in condition (C5) of \citet{xue.liang:2010} and condition (A5)
in \citet{ma.yang:2011}.

\subsection{Preliminary lemmas}\label{appa2}

Define $\alpha_{\ell}^{0}(\mathbf{x})=\sum_{k=1}^{d}\alpha
_{\ell
k}^{0}(x_{k})=B(\mathbf{x})^{\mathrm{T}}{\bolds{\gamma}}_{\ell}$,\break where $\alpha_{\ell k}^{0}(x_{k})$ is defined in (\ref
{DEF:glzl}%
). Let ${\bolds{\gamma}}_{I_{1}}=({\bolds{\gamma}}_{\ell}:\ell\in I_{1})^{\mathrm{T}}$. To prove Theorem~\ref
{THM:weakconsistencyLASSORE}, we next define the oracle estimator of $%
{\bolds{\gamma}}_{I_{1}}$ by minimizing the penalized
negative quasi-likelihood with all irrelevant predictors eliminated as such
%
\begin{equation}
L_{n}({\bolds{\gamma}}_{I_{1}})=\sum%
_{i=1}^{n}Q
\biggl[g^{-1}\biggl\{\sum_{\ell\in I_{1}}B(
\mathbf{X}%
_{i})^{\mathrm{T}}{\bolds{\gamma}}_{\ell}T_{\ell
}
\biggr\},Y_{i}\biggr]+n\lambda_{n}\sum
_{\ell\in I_{1}}w_{n\ell}\Vert {\bolds{\gamma}}_{\ell}
\Vert_{2}, \label{EQ:LASSORE}
\end{equation}
so that $\widehat{{\bolds{\gamma}}}_{I_{1}}^{0}=
({\widehat{{\bolds{\gamma}}}}_{\ell}^{0}:\ell\in
I_{1})^{\mathrm{T}%
}=\arg\min_{{\bolds{\gamma}}_{I_{1}}}L_{n}
({\bolds{\gamma}}_{I_{1}})$.
Define ${\widehat{{\bolds{\gamma}}}}_{I_{2}}^{0}=({\widehat{
{\bolds{\gamma}}}}_{\ell}^{0}:\ell\in I_{2})^{\mathrm{T}}$
with ${\widehat{{\bolds{\gamma}}}}_{\ell}^{0}\equiv
\mathbf
{0}%
_{dJ_{n}+1}$ for $\ell\in I_{2}$, where $\mathbf{0}_{dJ_{n}+1}$ is a $
(dJ_{n}+1)$-dimensional zero vector. We next present several lemmas, whose
detailed proofs are given in the online supplementary materials [\citet
{ma.carroll.liang.xu2015}]. Lemma~\ref%
{LEM:gammatilda} is used for the proof of Theorem~\ref%
{THM:weakconsistencyLASSORE}, while Lemma~\ref{LEM:Lngamma1} is needed
in the proof of Theorem~\ref{THM:oracle}.

\renewcommand{\thelemmaa}{\Alph{section}.\arabic{lemmaa}}
\begin{lemmaa}
\label{LEM:gammatilda}Under the conditions of Theorem~\ref%
{THM:weakconsistencyLASSORE}, one has
%
\begin{equation}
\bigl\Vert{\widehat{{\bolds{\gamma}}}}_{I_{1}}^{0}-{\bolds{
\gamma}}_{I_{1}}\bigr\Vert_{2}=O_{p}\bigl(
\lambda_{n}\Vert w_{n,I_{1}}\Vert+n^{-1/2}J_{n}^{1/2}+J_{n}^{-r}
\bigr), \label{RES1}
\end{equation}
and as $n\rightarrow\infty$,
%
\begin{equation}
P\bigl\{{\widehat{{\bolds{\gamma}}}=\bigl(\widehat{{\bolds{\gamma}}}}_{I_{1}}^{0\mathrm{T}},
{\widehat{{\bolds{\gamma}}}}%
_{I_{2}}^{0\mathrm{T}}
\bigr)^{\mathrm{T}}\bigr\}\rightarrow1. \label{RES2}
\end{equation}
\end{lemmaa}

\begin{lemmaa}
\label{LEM:Lngamma1}Under conditions \textup{(C1)--(C5)} and Assumptions \ref%
{Assumption1}--\ref{Assumption3},
%
\begin{equation}
\bigl\Vert{\widehat{{\bolds{\gamma}}}}_{,1}^{\mathcal
{S}}-%
{
\widehat{{\bolds{\gamma}}}}_{,1}^{\mathrm{OR}}\bigr\Vert
_{\infty
}=O_{p}\Bigl(\sqrt{\log n/\bigl(J_{n}^{\mathcal{S}}n
\bigr)}+\bigl(J_{n}^{\mathcal{S}%
}\bigr)^{-1/2}
\bigl(J_{n}^{\mathrm{ini}}\bigr)^{-r}\Bigr). \label{gammas-gammaor}
\end{equation}
\end{lemmaa}

\subsection{Proof of Theorem \texorpdfstring{\protect\ref{THM:weakconsistencyLASSORE}}{1}}\label{appa3}

By (\ref{RES1}) and (\ref{RES2}),
\begin{eqnarray*}
&&\sum_{\ell\in I_{1}}\Vert\widehat{\alpha}_{\ell}-
\alpha _{\ell}\Vert\asymp\Vert{\widehat{{\bolds{\gamma}}}}%
_{I_{1}}-{
\bolds{\gamma}}_{I_{1}}\Vert _{2}=O_{p}\bigl(\lambda
_{n}\Vert w_{n,I_{1}}\Vert+n^{-1/2}J_{n}^{1/2}+J_{n}^{-r}
\bigr),
\\
&&P\bigl(\llVert \widehat{\alpha}_{\ell}\rrVert >0,\ell\in
I_{1}\mbox{ and }\llVert \widehat{\alpha}_{\ell}\rrVert =0,
\ell\in I_{2}\bigr)\rightarrow1.
\end{eqnarray*}

\subsection{Proof of Theorem \texorpdfstring{\protect\ref{THM:asymptotics}}{2}}\label{appa4}
Let ${\bolds{\gamma}}_{,1}=
({\bolds{\gamma}}_{\ell1},\ell\in\widehat{I}_{1})^{\mathrm{T}}$, where
${\bolds{\gamma}}_{\ell1}$ is defined in~(\ref{EQ:spline}). By
Taylor's expansion, from (\ref{EQ:LnOR}), one has
\begin{eqnarray*}
{\widehat{{\bolds{\gamma}}}}_{,1}^{\mathrm
{OR}}-{\bolds{
\gamma}}_{,1} &=&\Biggl[\sum_{i=1}^{n}Z_{i,1}Z_{i,1}^{%
\mathrm{T}}
\bigl\{\dot{g}^{-1}\bigl(\eta_{i}^{\ast}\bigr)\bigr
\}^{2}\Big/\sigma_{i}^{2}\Biggr]^{-1}
\\
&&{}\times\Biggl[\sum_{i=1}^{n}Z_{i,1}
\bigl\{ Y_{i}-g^{-1}\bigl(\eta _{i}^{0}
\bigr)\bigr\}\bigl(\dot{g}^{-1}\bigl(\eta_{i}^{0}
\bigr)/\sigma_{i}^{2}\bigr)\Biggr],
\end{eqnarray*}
where $\eta_{i}^{0}=\sum_{\ell=1}^{p}\{\alpha_{\ell
0}+\sum_{k=2}^{d}\alpha_{\ell k}(X_{ik})\}T_{i\ell
}+\sum_{\ell=1}^{p}B^{\mathcal{S}}(x_{1})^{\mathrm{T}}%
{\bolds{\gamma}}_{\ell1}T_{i\ell}$ and
\[
\eta_{i}^{\ast}=\sum_{\ell=1}^{p}
\Biggl\{\alpha_{\ell
0}+\sum_{k=2}^{d}
\alpha_{\ell k}(X_{ik})\Biggr\}T_{i\ell
}+\sum
_{\ell=1}^{p}B^{\mathcal{S}}(x_{1})^{\mathrm{T}}%
{
\bolds{\gamma}}_{\ell1}^{\ast}T_{i\ell},
\]
where ${\bolds{\gamma}}_{,1}^{\ast}=(%
{\bolds{\gamma}}_{\ell1}^{\ast},\ell\in\widehat{I}_{1})^{\mathrm{T}}\in
({\bolds{\gamma}}_{,1},{\widehat{{\bolds{\gamma}}}}_{,1}^{\mathrm{OR}})$.
Following similar reasoning as the proofs
for (\ref{RES1}), we have
$\Vert{\widehat{{\bolds{\gamma}}}}_{,1}^{\mathrm{OR}}-
{\bolds{\gamma}}_{,1}\Vert_{2}=o_{p}(1)$. Then
${\widehat{{\bolds{\gamma}}}}_{,1}^{\mathrm{OR}}
-{\bolds{\gamma}}_{,1}= (
{\widehat{{\bolds{\gamma}}}}_{,1e}^{\mathrm
{OR}}+{\widehat{{\bolds{\gamma}}}}_{,1\mu}^{\mathrm{OR}} )
+o_{p}(1)$,
where
%
\begin{eqnarray}
\label{EQ:zzd} {\widehat{{\bolds{\gamma}}}}_{,1e}^{\mathrm{OR}} &=&
\Biggl[\sum_{i=1}^{n}Z_{i,1}Z_{i,1}^{\mathrm{T}}
\bigl\{\dot{g}^{-1}(\eta _{i})\bigr\}^{2}\Big/
\sigma_{i}^{2}\Biggr]^{-1}\Biggl[\sum
_{i=1}^{n}Z_{i,1}\varepsilon
_{i}\bigl\{\dot{g}^{-1}(\eta_{i})/
\sigma_{i}^{2}\bigr\}\Biggr],
\nonumber
\\
{\widehat{{\bolds{\gamma}}}}_{,1\mu}^{\mathrm{OR}} &=&\Biggl[\sum
_{i=1}^{n}Z_{i,1}Z_{i,1}^{\mathrm{T}}
\bigl\{\dot{g}^{-1}(\eta _{i})\bigr\}^{2}\Big/
\sigma_{i}^{2}\Biggr]^{-1}
\\
&&\times{}\Biggl[\sum_{i=1}^{n}Z_{i,1}
\bigl\{g^{-1}(\eta _{i})-g^{-1}\bigl(
\eta_{i}^{0}\bigr)\bigr\} \bigl\{\dot{g}^{-1}(
\eta_{i})/\sigma _{i}^{2}\bigr\}\Biggr].\nonumber
\end{eqnarray}
Therefore, $\operatorname{var}(\widehat{\bolds{\gamma}}_{,1e}^{\mathrm{OR}%
}|\mathbb{X},\mathbb{T})=[\sum_{i=1}^{n}Z_{i,1}Z_{i,1}^{\mathrm
{T}}\{%
\dot{g}^{-1}(\eta_{i})\}^{2}/\sigma_{i}^{2}]^{-1}$. By Theorem~5.4.2 of
\citet{DeVore.Lorentz:1993}, for sufficiently large $n$, there exist
constants $0<c_{B}\leq C_{B}<\infty$, such that $c_{B}\mathbf
{I}_{J_{n}^{%
\mathcal{S}}\times J_{n}^{\mathcal{S}}}\leq E(B_{1}^{\mathcal{S}%
}(X_{i1})B_{1}^{\mathcal{S}}(X_{i1})^{\mathrm{T}})\leq C_{B}\mathbf
{I}_{J_{n}^{%
\mathcal{S}}\times J_{n}^{\mathcal{S}}}$. By condition (C5), for $n$ large
enough, there are constants $0<C_{T},C^{\prime}<\infty$, such
that
\begin{eqnarray*}
&&E\bigl[Z_{i,1}Z_{i,1}^{\mathrm{T}}\bigl\{
\dot{g}^{-1}(\eta_{i})\bigr\}^{2}/\sigma
_{i}^{2}\bigr]\\
&&\qquad\leq C^{\prime}E\bigl[\bigl
\{B_{1}^{\mathcal
{S}}(X_{i1})B_{1}^{\mathcal
{S}%
}(X_{i1})^{\mathrm{T}}
\bigr\}\otimes\bigl\{E(T_{\ell}T_{\ell^{\prime}}|\mathbf {X}%
)
\bigr\}_{\ell,\ell^{\prime}\in\widehat{I}_{1}}\bigr]
\\
&&\qquad \leq CC_{T}s^{\ast}E\bigl\{B_{1}^{\mathcal{S}}(X_{i1})B_{1}^{\mathcal{%
S}}(X_{i1})^{\mathrm{T}}
\bigr\}\otimes\mathbf{I}_{s^{\ast}\times s^{\ast
}}\leq C^{\prime}C_{T}C_{B}s^{\ast}
\mathbf{I}_{J_{n}^{\mathcal{S}}\times
J_{n}^{%
\mathcal{S}}}\otimes\mathbf{I}_{s^{\ast}\times s^{\ast}}
\\
&&\qquad =Cs^{\ast}%
\mathbf{I}_{J_{n}^{\mathcal{S}}s^{\ast}\times J_{n}^{\mathcal
{S}}s^{\ast}},
\end{eqnarray*}
where $C=$ $C^{\prime}C_{T}C_{B}$. Similarly, we have $%
E[Z_{i,1}Z_{i,1}^{\mathrm{T}}\{\dot{g}^{-1}(\eta_{i})\}^{2}/\sigma
_{i}^{2}]\geq \break  c\mathbf{I}_{J_{n}^{\mathcal{S}}s^{\ast}\times J_{n}^{%
\mathcal{S}}s^{\ast}}$ for some constant\vspace*{1pt} $0<c<\infty$. Thus,
following the
same reasoning as the proof for (S.5) in the supplementary
materials [\citet{ma.carroll.liang.xu2015}], we have
with probability $1$, for $n\rightarrow\infty$,
%
\begin{eqnarray}
\label{EQ:ZZdsig}
C^{-1}\bigl(s^{\ast}
\bigr)^{-1}n^{-1}\mathbf{I}_{J_{n}^{\mathcal{S}}s^{\ast
}\times
J_{n}^{\mathcal{S}}s^{\ast}}&\leq& {}\Biggl[ \sum
_{i=1}^{n}Z_{i,1}Z_{i,1}^{\mathrm{T}}
\bigl\{\dot{g}^{-1}(\eta _{i})\bigr\}^{2}\Big/
\sigma_{i}^{2}\Biggr]^{-1}
\nonumber
\\[-8pt]
\\[-8pt]
\nonumber
&\leq& c^{-1}n^{-1}\mathbf{I}_{J_{n}^{%
\mathcal{S}}s^{\ast}\times J_{n}^{\mathcal{S}}s^{\ast}}.
\end{eqnarray}
%
By the Lindeberg central limit theorem, it can be proved that
%
\begin{equation}
\mathbf{a}^{\mathrm{T}}\sigma_{n}^{-1}(x_{1})
\bigl\{B^{\mathcal
{S}}(x_{1})%
\widehat{{\bolds{
\gamma}}}_{,1e}^{\mathrm{OR}}\bigr\}\rightarrow N(0,1), \label{EQ:nor}
\end{equation}
for any $\mathbf{a}\in R^{s^{\ast}}$ with $\Vert\mathbf{a}\Vert_{2}=1$.
Since $\mathbf{a}^{\mathrm{T}}\sigma_{n}^{-1}(x_{1})\{\widehat{\alpha
}_{1}^{%
\mathrm{OR}}(x_{1})-b_{1}(x_{1})\}= \mathbf{a}^{\mathrm{T}}\times\break \sigma
_{n}^{-1}(x_{1})\{B^{\mathcal{S}}(x_{1})
{\widehat{{\bolds{\gamma}}}}_{,1e}^{\mathrm{OR}}\}+o_{p} ( 1 ) $, by (\ref{EQ:nor})
and Slutsky's theorem, we have
%
\begin{equation}
\mathbf{a}^{\mathrm{T}}\sigma_{n}^{-1}(x_{1})
\bigl\{\widehat{\alpha }_{1}^{\mathrm{OR%
}}(x_{1})-b_{1}(x_{1})
\bigr\}\rightarrow N(0,1). \label{EQ:normality}
\end{equation}
By (\ref{EQ:zzd}) and (\ref{EQ:ZZdsig}), with probability approaching
$1$,
\begin{eqnarray*}
\sum_{\ell\in I_{1}}\bigl\Vert\widehat{\alpha}_{\ell1}^{\mathrm
{OR}%
}-b_{\ell1}
\bigl\Vert^{2} &\asymp&\bigl\Vert{\widehat{{\bolds{\gamma}}}}_{,1e}^{\mathrm{OR}}
\bigr\Vert_{2}^{2}
\\[-2pt]
& \leq& c^{-2}n^{-2}\Biggl[\sum
_{i=1}^{n}\varepsilon _{i}Z_{i,1}^{\mathrm{T}}
\bigl(\dot{g}^{-1}(\eta_{i})/\sigma _{i}^{2}
\bigr)\Biggr] \Biggl[\sum_{i=1}^{n}Z_{i,1}
\varepsilon_{i}\bigl(\dot {g}^{-1}(\eta _{i})/
\sigma_{i}^{2}\bigr)\Biggr]
\\
&\asymp& c^{-2}n^{-1}E\bigl[Z_{i,1}^{\mathrm{T}}Z_{i,1}
\bigl\{\dot {g}^{-1}(\eta _{i})\bigr\}^{2}/
\sigma_{i}^{2}\bigr]\asymp s^{\ast}J_{n}^{\mathcal{S}}n^{-1};
\\
\bigl\Vert\mathbf{a}^{\mathrm{T}}\bigl(\widehat{\alpha}_{1}^{\mathrm
{OR}}-b_{1}
\bigr)\bigr\Vert ^{2} &\leq&C_{a}\bigl\Vert{\widehat{{\bolds{
\gamma}}}}_{,1e}^{%
\mathrm{OR}}\bigr\Vert_{2}^{2}\leq
C_{a}c^{-1}n^{-2}\Biggl(\sum
_{i=1}^{n}\varepsilon_{i}Z_{i,1}^{\mathrm
{T}%
}
\Biggr) \Biggl(\sum_{i=1}^{n}Z_{i,1}
\varepsilon_{i}\Biggr)
\\
&\asymp&C_{a}c^{-1}n^{-1}E\bigl(Z_{i,1}^{\mathrm{T}}Z_{i,1}
\bigr)\asymp s^{\ast
}J_{n}^{%
\mathcal{S}}n^{-1}.
\end{eqnarray*}
Since $\sup_{x_{1}\in{}[0,1]}|\alpha_{\ell
1}(x_{1})-B_{1}^{\mathcal{S%
}}(x_{1})^{\mathrm{T}}{\bolds{\gamma}}_{\ell1}|=O\{
(J_{n}^{\mathcal{S%
}})^{-r}\}$, it can be proved that
$\Vert\mathbf{a}^{\mathrm{T}}{\widehat{{\bolds{\gamma}}}}
_{,1\mu}^{\mathrm{OR}}\Vert\leq\Vert{\widehat{{\bolds{\gamma}}}}_{,1\mu}^{\mathrm{OR}}\Vert_{2}=O_{p}\{(s^{\ast
})^{1/2}(J_{n}^{%
\mathcal{S}})^{-r}\}$, and $\Vert\mathbf{a}^{\mathrm{T}}(b_{1}-\alpha
_{1}^{0})\Vert\asymp
\Vert\mathbf{a}^{\mathrm{T}}{\widehat{{\bolds{\gamma}}}}%
_{,1\mu2}^{\mathrm{OR}}\Vert=O_{p}\{(s^{\ast})^{1/2}(J_{n}^{\mathcal{S}
})^{-r}\}$. Hence
\[
\bigl\Vert\mathbf{a}^{\mathrm{T}}(b_{1}-\alpha_{1})\bigr\Vert\leq
\bigl\Vert \mathbf {a}^{%
\mathrm{T}}\bigl(b_{1}-\alpha_{1}^{0}
\bigr)\bigr\Vert+\bigl\Vert\mathbf{a}^{\mathrm
{T}}\bigl(\alpha _{1}^{0}-
\alpha_{1}\bigr)\bigr\Vert=O_{p}\bigl\{s^{\ast}
\bigl(J_{n}^{\mathcal
{S}}\bigr)^{-r}\bigr\}.
\]
By (\ref{EQ:normality}), $\{\mathbf{e}_{\ell}^{\mathrm{T}}\sigma
_{n}^{2}(x_{1})\mathbf{e}_{\ell}\}^{-1/2}\{\widehat{\alpha}_{\ell1}^{
\mathrm{OR}}(x_{1})-b_{\ell1}\}\rightarrow N(0,1)$, and $\sup_{\ell
\in
\widehat{I}_{1}}|\widehat{\alpha}_{\ell0}^{\mathrm{OR}}-\alpha
_{\ell
0}|=O_{p}(n^{-1/2})$ follows from the central limit theorem.

\subsection{Proof of Theorem \texorpdfstring{\protect\ref{THM:oracle}}{3}}\label{appa5}

By (\ref{gammas-gammaor}) in Lemma~\ref{LEM:Lngamma1},
\[
\sup_{x_{1}\in{}[0,1]}\bigl\Vert\widehat{\alpha }_{1}^{\mathcal{S}%
}(x_{1})-
\widehat{\alpha}_{1}^{\mathrm{OR}}(x_{1})
\bigr\Vert_{\infty}\leq \sup_{x_{1}\in{}[0,1]}\sum
_{j=1}^{J_{n}^{\mathcal{S}%
}}\bigl|B_{j,1}^{\mathcal{S}}(x_{1})\bigr|
\bigl\Vert{\widehat{{\bolds{\gamma}}}}_{,1}^{\mathcal{S}}- {\widehat{{
\bolds{\gamma}}}}%
_{,1}^{\mathrm{OR}}\bigr\Vert_{\infty}.
\]
The right-hand side is bounded by $O_{p}\{(n^{-1}\log
n)^{1/2}+(J_{n}^{\mathrm{ini}})^{-r}\}$.
$\Vert\widehat{\bolds{\alpha}}_{0}^{\mathcal
{S}}-\widehat{\bolds{\alpha}}_{0}^{\mathrm{OR}}\Vert_{2}=o_{p}(n^{-1/2})$
can be
proved following the same procedure and thus omitted. By (\ref{EQ:ZZdsig}),
with probability approaching $1$, for large enough $n$, for any
$x_{1}\in
{}[0,1]$, and $\mathbf{a}\in R^{s^{\ast}}$ with $\Vert\mathbf
{a}\Vert
_{2}=1$, one has
\begin{eqnarray*}
\mathbf{a}^{\mathrm{T}}\sigma_{n}^{2}(x_{1})
\mathbf{a} &\leq &c_{Z}^{-1}n^{-1}%
\mathbf{a}^{\mathrm{T}}B^{\mathcal{S}}(x_{1})B^{\mathcal
{S}}(x_{1})^{\mathrm
{T}}%
\mathbf{a}\leq c^{-1}J_{n}^{\mathcal{S}}n^{-1}
\mathbf{a}^{\mathrm
{T}}\mathbf{a}%
,
\\
\mathbf{a}^{\mathrm{T}}\sigma_{n}^{2}(x_{1})
\mathbf{a} &\geq &C_{Z}^{-1}\bigl(s^{\ast}
\bigr)^{-1}n^{-1}\mathbf{a}^{\mathrm{T}}B^{\mathcal{S}%
}(x_{1})B^{\mathcal{S}}(x_{1})^{\mathrm{T}}
\mathbf{a}\geq C^{-1}J_{n}^{
\mathcal{S}}\bigl(s^{\ast}
\bigr)^{-1}n^{-1}\mathbf{a}^{\mathrm{T}}\mathbf{a},
\end{eqnarray*}
where $\sigma_{n}^{2}(x_{1})$ is defined in (\ref{DEF:sig2nx1}). Thus
\begin{eqnarray*}
&&\sup_{x_{1}\in{}[0,1]}\bigl|\mathbf{a}^{\mathrm{T}}\sigma
_{n}^{-1}(x_{1})\bigl\{\widehat{
\alpha}_{1}^{\mathcal{S}}(x_{1})-\widehat {\alpha
}_{1}^{\mathrm{OR}}(x_{1})\bigr\}\bigr|
\\
&&\qquad\leq\sup_{x_{1}\in{}[0,1]}\bigl\Vert\sigma _{n}^{-1}(x_{1})
\bigr\Vert_{2}\bigl\Vert\widehat{\alpha}_{1}^{\mathcal
{S}}(x_{1})-%
\widehat{\alpha}_{1}^{\mathrm{OR}}(x_{1})
\bigr\Vert_{2}
\\
&&\qquad=O_{p}\bigl[s^{\ast}\bigl\{\bigl(\log
n/J_{n}^{\mathcal
{S}}\bigr)^{1/2}+\bigl(n/J_{n}^{%
\mathcal{S}}
\bigr)^{1/2}\bigl(J_{n}^{\mathrm{ini}}\bigr)^{-r}
\bigr\}\bigr]=o_{p}(1).
\end{eqnarray*}

\subsection{Proof of Theorem \texorpdfstring{\protect\ref{THM:bands}}{4}}

Using the strong approximation lemma given in Theorem~2.6.7 of
\citet{Csorgo.Revesz:1981}, we can prove by the same procedure as
Lemma A.7
in \citet{ma:yang:carroll:2012} that
%
\begin{equation}
\sup_{x_{1}\in{}[0,1]}\bigl\llvert \widehat{\alpha}_{\ell
1}^{%
\mathrm{OR}}(x_{1})-b_{\ell1}(x_{1})-
\widehat{\alpha}_{\ell
1,\varepsilon
}^{0}(x_{1})\bigr\rrvert
=o_{\mathrm{a.s.}}\bigl(n^{t}\bigr) \label{EQ:mhate1}
\end{equation}
for some $t<-r/(2r+1)<0$, where $\widehat{\alpha}_{\ell1,\varepsilon
}^{0}(x_{1})$ is
\[
\mathbf{e}_{\ell}^{\mathrm{%
T}}B^{\mathcal{S}}(x_{1})
\Biggl[\sum_{i=1}^{n}Z_{i,1}Z_{i,1}^{\mathrm
{T}}
\bigl\{%
\dot{g}^{-1}(\eta_{i})\bigr
\}^{2}\Big/\sigma _{i}^{2}\Biggr]^{-1}\Biggl[
\sum_{i=1}^{n}Z_{i,1}e_{i}
\bigl\{\dot{g}^{-1}(\eta _{i})/\sigma_{i}^{2}
\bigr\}\Biggr],
\]
and $e_{i},1\leq i\leq n$, are i.i.d. $N(0,1)$ independent of $
Z_{i,1}$. For $\sigma_{n}^{2}(x_{1})$ defined in~(\ref{DEF:sig2nx1})
and $%
\sigma_{n1}(x_{1})\asymp(J_{n}^{\mathcal{S}}/n)^{1/2}\{1+o_{p}(1)\}$
uniformly in $x_{1}\in{}[0,1]$. By (\ref{EQ:mhate1}),
$J_{n}^{\mathcal{%
S}}\asymp n^{1/(2r+1)}$ and $t<-r/(2r+1)<0$, we have
%
\begin{eqnarray}
&&\sup_{x_{1}\in[\label{EQ:mhate}0,1]}\bigl\vert\bigl\{\log(L_{n}+1)\bigr
\}^{-1/2}\sigma _{n1}^{-1}(x_{1})\bigl\{
\widehat{\alpha}_{\ell1}^{\mathrm
{OR}}(x_{1})-b_{\ell
1}(x_{1})-
\widehat{\alpha}_{\ell1,\varepsilon}^{0}(x_{1})\bigr\}\bigr\vert
\nonumber
\\
&&\qquad=o_{\mathrm{a.s.}}\bigl(\bigl\{\log(L_{n}+1)\bigr\}^{-1/2}
\bigl(n/J_{n}^{\mathcal{S}%
}\bigr)^{1/2}n^{t}
\bigr)
\\
&&\qquad=o_{\mathrm{a.s.}}\bigl(\bigl\{\log(L_{n}+1)\bigr
\}^{-1/2}n^{r/(2r+1)-t}\bigr)=o_{\mathrm{a.s.}}(1).\nonumber
\end{eqnarray}
Define $\eta(x_{1})=\sigma_{n1}^{-1}(x_{1})\widehat{\alpha}_{\ell
1,\varepsilon}^{0}(x_{1})$. It is apparent that $\mathcal{L\{}\eta
(\xi
_{J})\vert Z_{i,1},1\leq i\leq n \}=N(0,1)$, so $\mathcal
{L\{}%
\eta(\xi_{J})\}=N(0,1)$ for $0\leq J\leq L_{n}$. Moreover, the eigenvalues
of $(EZ_{i,1}Z_{i,1}^{\mathrm{T}})^{-1}\asymp J_{n}^{\mathcal{S}}$.
Then with
probability approaching $1$, for $J\neq J^{\prime}$,
\begin{eqnarray*}
\bigl\llvert E\mathcal{\bigl\{}\eta(\xi_{J})\eta(
\xi_{J^{\prime
}})\bigr\}\bigr\rrvert &\asymp&\bigl(n/J_{n}^{\mathcal{S}}
\bigr)n^{-1}\bigl\llvert \mathbf {e}%
_{\ell}^{\mathrm{T}}B^{\mathcal{S}}(
\xi_{J}) \bigl(EZ_{i,1}Z_{i,1}^{\mathrm
{T}%
}
\bigr)^{-1}B^{\mathcal{S}}(\xi_{J^{\prime}})^{\mathrm{T}}\mathbf
{e}_{\ell
}\bigr\rrvert
\\
&\asymp&\bigl\llvert \mathbf{e}_{\ell}^{\mathrm
{T}}B^{\mathcal{S}}(
\xi_{J})B^{%
\mathcal{S}}(\xi_{J^{\prime}})^{\mathrm{T}}
\mathbf{e}_{\ell}\bigr\rrvert =\sum_{j=1}^{J_{n}^{\mathcal{S}}}B_{j,1}^{\mathcal{S}}(
\xi _{J})B_{j,1}^{\mathcal{S}}(\xi_{J^{\prime}})
\end{eqnarray*}
and $\sum_{j=1}^{J_{n}^{\mathcal{S}}}B_{j,1}^{\mathcal{S}}(\xi
_{J})B_{j,1}^{%
\mathcal{S}}(\xi_{J^{\prime}})\asymp C$ for a constant $0<C<\infty$
when $%
\llvert  j_{J}-j_{J^{\prime}}\rrvert \leq(q-1)$ and $%
\sum_{j=1}^{J_{n}^{\mathcal{S}}}B_{j,1}^{\mathcal{S}}(\xi_{J})B_{j,1}^{
\mathcal{S}}(\xi_{J^{\prime}})=0$ when $\llvert
j_{J}-j_{J^{\prime
}}\rrvert >(q-1)$, in which $j_{J}$ denotes the index of the knot
closest to $\xi_{J}$ from the left. Therefore, by $L_{n}\asymp J_{n}^{%
\mathcal{S}}$, there exist constants $0<C_{1}<\infty$ and
$0<C_{2}<\infty$
such that with probability approaching $1$, for $J\neq J^{\prime}$, $%
\llvert  E\mathcal{\{}\eta(\xi_{J})\eta(\xi_{J^{\prime
}})\}\rrvert \leq C_{1}^{-\llvert  j_{J}-j_{J^{\prime}}\rrvert
}\leq C_{2}^{-\llvert  J-J^{\prime}\rrvert }$. By Lemma~A1
given in
\citet{ma:yang:2011}, we have
\[
\lim_{n\rightarrow\infty}P\Bigl\{\sup_{0\leq J\leq
L_{n}}\bigl\vert\bigl\{2
\log(L_{n}+1)\bigr\}^{-1/2}\eta(\xi_{J})\bigr\vert\leq
d_{N_{n}}(\alpha)\Bigr\}=1-\alpha,
\]
and hence
%
\begin{equation}\qquad
\lim_{n\rightarrow\infty}P\Bigl\{\sup_{x_{1}\in S_{n,1}}\bigl\llvert
\bigl\{2\log(L_{n}+1)\bigr\}^{-1/2}\sigma_{n1}^{-1}(x_{1})
\widehat{\alpha }_{\ell
1,\varepsilon}^{0}(x_{1})\bigr\rrvert
\leq d_{N_{n}}(\alpha)\Bigr\} =1-\alpha. \label{EQ:extreme}
\end{equation}
Furthermore, according to the result on page 149 of \citet{deboor:2001},
we have
%
\begin{eqnarray}
\label{EQ:m1ORm} &&\sup_{x_{1}\in{}[0,1]}\bigl\llvert \bigl\{\log
(L_{n}+1)\bigr\}^{-1/2}\sigma_{n1}^{-1}(x_{1})
\bigl\{b_{\ell1}(x_{1})-\alpha _{\ell
1}(x_{1})
\bigr\}\bigr\rrvert
\nonumber
\\[-8pt]
\\[-8pt]
\nonumber
&&\qquad=O_{p}\bigl(\bigl\{\log(L_{n}+1)\bigr\}^{-1/2}
\bigl(n/J_{n}^{\mathcal
{S}}\bigr)^{1/2}\bigl(J_{n}^{%
\mathcal{S}}
\bigr)^{-r}\bigr)=o_{p}(1).
\end{eqnarray}
Moreover,
$
\widehat{\alpha}_{\ell1}^{\mathrm{OR}}(x_{1})-\alpha_{\ell1}(x_{1})=
\widehat{\alpha}_{\ell1,\varepsilon}^{0}(x_{1})+\{\widehat{\alpha
}_{\ell
1}^{\mathrm{OR}}(x_{1})-b_{\ell1}(x_{1})-\widehat{\alpha}_{\ell
1,\varepsilon}^{0}(x_{1})\}+\{b_{\ell1}(x_{1})-\alpha_{\ell
1}(x_{1})\}$.
Hence by (\ref{EQ:mhate}) and (\ref{EQ:m1ORm}), we have
%
\begin{eqnarray}
\label{EQ:supalpha0} && \lim_{n\rightarrow\infty}P \Bigl\{ \sup
_{x_{1}\in
S_{n,1}}\bigl\{\log (L_{n}+1)\bigr\}^{-1/2}
\sigma_{n1}^{-1}(x_{1})\bigl\llvert \widehat{
\alpha }_{\ell
1}^{\mathrm{OR}}(x_{1})-\alpha_{\ell1}(x_{1})
\bigr\rrvert \leq d_{N_{n}}(\alpha) \Bigr\}
\nonumber
\\
&&\qquad =\lim_{n\rightarrow\infty}P \Bigl\{ \sup_{x_{1}\in
S_{n,1}}\bigl\{
\log(L_{n}+1)\bigr\}^{-1/2}\sigma_{n1}^{-1}(x_{1})
\bigl\llvert \widehat{%
\alpha}_{\ell1,\varepsilon}^{0}(x_{1})
\bigr\rrvert \leq d_{N_{n}}(\alpha ) \Bigr\}
\\
&&\qquad =1-\alpha,\nonumber
\end{eqnarray}
where the last step follows from \eqref{EQ:extreme}. By the oracle property
given in Theorem~\ref{THM:oracle}, and $J_{n}^{\mathcal{S}}\asymp
n^{1/(2r+1)}$ and $n^{1/(2r+1)}\ll J_{n}^{\mathrm{ini}}$, we have
%
\begin{eqnarray}\qquad
\label{EQ:alphas} &&\sup_{x_{1}\in{}[0,1]}\bigl\{\log(L_{n}+1)\bigr
\}^{-1/2}\sigma _{n1}^{-1}(x_{1})\bigl|
\widehat{\alpha}_{\ell1}^{\mathcal
{S}}(x_{1})-
\widehat{%
\alpha}_{\ell1}^{\mathrm{OR}}(x_{1})\bigr|
\nonumber
\\[-8pt]
\\[-8pt]
\nonumber
&&\qquad =O_{p}\bigl[\log(L_{n}+1)^{-1/2}
\bigl(n/J_{n}^{\mathcal
{S}}\bigr)^{1/2}\bigl(n^{-1}
\log n\bigr)^{1/2}+\bigl(J_{n}^{\mathrm{ini}}
\bigr)^{-r}\bigr]=o_{p}(1).
\end{eqnarray}
Therefore, by (\ref{EQ:supalpha0}) and (\ref{EQ:alphas}), we have
\[
\lim_{n\rightarrow\infty}P \Bigl\{ \sup_{x_{1}\in S_{n,1}}\bigl\{ \log
(L_{n}+1)\bigr\}^{-1/2}\sigma_{n1}^{-1}(x_{1})
\bigl\llvert \widehat{\alpha }_{\ell
1}^{\mathcal{S}}(x_{1})-
\alpha_{\ell1}(x_{1})\bigr\rrvert \leq d_{N_{n}}(
\alpha) \Bigr\} =1-\alpha,
\]
and hence the result in Theorem~\ref{THM:bands} is proved.
\end{appendix}

\section*{Acknowledgments}
The authors thank the Co-Editors, an Associate Editor and three referees
for their valuable suggestions and comments that have substantially
improved an earlier version of this paper.

\begin{supplement}[id=suppA]
\stitle{Supplemental materials for ``Estimation and inference in generalized
additive coefficient models for nonlinear interactions with high-dimensional
covariates''}
\slink[doi]{10.1214/15-AOS1344SUPP} 
\sdatatype{.pdf}
\sfilename{aos1344\_supp.pdf}
\sdescription{The supplementary material presents additional numerical
results and the proofs of Lemmas \ref{LEM:gammatilda} and \ref{LEM:Lngamma1}.}
\end{supplement}

%




\printaddresses
\end{document}